\documentclass[10pt]{article}
\usepackage{eufrak}
\usepackage{euscript}
\usepackage{pictex}

\newtheorem{lemma}{{\sc LEMMA}}[section]
\newtheorem{theorem}{{\sc THEOREM}}[section]

\newtheorem{assumption}{Assumption}[section]
\newtheorem{procedure}{Discretization procedure}[section]

\newtheorem{corollary}{{\sc COROLLARY}}[section]
\newtheorem{Figure}{Figure}

\newcommand {\Proof}  {\mbox{\sc Proof: }}
\newcommand {\QED}{\hfill {\bf \textsf{QED}}\\}

\def\1{\hbox{\upshape1\kern-.15em\vrule height 1.6ex width .3pt
\vrule width .8pt height .25pt\kern.15em}}
\def\lev{\langle \,}
\def\des{\,\rangle}
\def\lnorm{\,\rule[-1mm]{0.6mm}{4mm}\,}
\def\blnorm{ \,\rule[-2mm]{0.6mm}{7mm}\,}
\def\der{\,\rule[0mm]{0.1mm}{1.3mm} \rule[0mm]{1.3mm}{0.1mm}
 \hspace{-1.3mm} \rule[1.3mm]{1.3mm}{0.1mm}
 \hspace{-1.3mm} \rule[2.3mm]{1.3mm}{0.1mm}
 \rule[0mm]{0.1mm}{2.3mm} \,}

\def\<<{\,<\hspace{-1.3mm}<}
\def\>>{>\hspace{-1.3mm}>}
\def\notsubseteq{\subseteq \hspace{-3.3mm}/\hspace{2.0mm}}

\def\msb{\mbox{\boldmath{$b$}}}
\def\msx{\mbox{\boldmath{$x$}}}
\def\msxd{\mbox{\boldmath{\scriptsize$x$}}}
\def\msy{\mbox{\boldmath{$y$}}}
\def\msyd{\mbox{\boldmath{\scriptsize$y$}}}
\def\msv{\mbox{\boldmath{$v$}}}

\def\msw{\mbox{\boldmath{$w$}}}
\def\mswd{\mbox{\boldmath{\scriptsize$w$}}}
\def\msz{\mbox{\boldmath{$z$}}}
\def\mszd{\mbox{\boldmath{\scriptsize$z$}}}
\def\msp{\mbox{\boldmath{$p$}}}

\def\msq{\mbox{\boldmath{$q$}}}
\def\msr{\mbox{\boldmath{$r$}}}
\def\msrd{\mbox{\boldmath{\scriptsize$r$}}}
\def\msk{\mbox{\boldmath{$k$}}}
\def\mskd{\mbox{\boldmath{\scriptsize$k$}}}
\def\msm{\mbox{\boldmath{$m$}}}
\def\msmd{\mbox{\boldmath{\scriptsize$m$}}}
\def\msl{\mbox{\boldmath{$l$}}}
\def\msld{\mbox{\boldmath{\scriptsize$l$}}}
\def\mse{\mbox{\boldmath{$e$}}}
\def\msed{\mbox{\boldmath{\scriptsize$e$}}}

\def\mst{\mbox{\boldmath{$t$}}}
\def\mss{\mbox{\boldmath{$s$}}}

\font\msbm = msbm10
\def\Msbm#1{\hbox{\msbm #1}}

\newcommand {\Mg} {\overline{M}}
\newcommand {\Md} {\underline{M}}
\newcommand {\bbJ}{\1}
\newcommand {\bbR}{{\Msbm R}}
\newcommand {\bbN}{{\Msbm N}}
\newcommand {\bbZ}{{\Msbm Z}}

\newdimen\vskp
\def\AMSclass#1{\vskip\vskp\vbox{\hbox{\vbox{\def\ams{#1}%
      \ifx\ams\empty
      \errmessage{Please put AMS subject classification in \noexpand
      \AMSclass command}
      \else\noindent
      \small\strut{\bf AMS subject classification:\hskip0.5em}\ams
      \strut}}}\fi}

\def\KeyWords#1{\vskip\vskp\vbox{\hbox{\vbox{\def\Kwd{#1}%
      \ifx\Kwd\empty
      \errmessage{Please fill Key words in \noexpand
      \KeyWords command}
      \else
      \newdimen\keywdswidth
      \setbox0=\hbox{\small\bf Key words:\hskip0.5em}
      \keywdswidth=\wd0
      \noindent\hangindent=\keywdswidth
      \hangafter=1
      \small\strut{\bf Key words:\hskip0.5em}\Kwd\strut}}}\fi}
\title{MARKOV JUMP PROCESSES APPROXIMATING A NON-SYMMETRIC GENERALIZED
DIFFUSION: NUMERICS EXPLAINED TO PROBABILISTS\thanks{Supported
      by grant 0037014 of the Ministry of Science, Higher Education and Sports,
Croatia.}}

\author{Ned\v zad Limi\'c
      \thanks{
      Dept. of Mathematics, University of Zagreb,
      Bijeni\v{c}ka 30, 10002 Zagreb, Croatia,
      e--mail: nlimic@math.hr}}

\begin{document}
\maketitle
\begin{abstract}
\noindent
Consider a non-symmetric generalized diffusion $X(\cdot)$ in ${\bbR}^d$
determined by the differential operator $A(\msx)=-\sum_{ij}
\partial_ia_{ij}(\msx)\partial_j +\sum_i b_i(\msx)\partial_i$. In this paper
the diffusion process is approximated by Markov jump processes $X_n(\cdot)$,
in homogeneous and isotropic grids $G_n \subset {\bbR}^d$, which converge in
distribution to the diffusion $X(\cdot)$. The generators of $X_n(\cdot)$ are
constructed explicitly. Due to the homogeneity and isotropy of grids, the proposed
method for $d\geq3$ can be applied to processes for which the diffusion tensor
$\{a_{ij}(\msx)\}_{11}^{dd}$ fulfills an additional condition.
The proposed construction offers a simple method for simulation of sample paths
of non-symmetric generalized diffusion. Simulations are carried out in terms of
jump processes $X_n(\cdot)$. For $d=2$ the construction can be easily
implemented into a computer code.
\end{abstract}
\AMSclass{60H35, 60J22, 60J27, 60J60, 65C30}

\section{INTRODUCTION}\label{sec1}
A symmetric tensor valued function $\msx \mapsto a(\msx)=\{a_{ij}(\msx)
\}_{11}^{dd}$ which is measurable, bounded and strictly positive definite on
${\bbR}^d$ defines a second order differential operator in divergence form
on ${\bbR}^d$, $A_0(\msx)=-\sum_{ij}\partial_i a_{ij}(\msx)\partial_j$. Each
$A_0(\msx)$ determines a symmetric diffusion $X(\cdot)$ in ${\bbR}^d$. In~\cite{SZ}
the process $X(\cdot)$ is approximated by a sequence of Markov jump processes
$X_n(\cdot)$ (MJP) in grids, i.e.~processes in continuous time and discrete state
spaces. Coordinate systems in in~\cite{SZ}
depend locally on the structure of the tensor valued
function $a$. The object of the present analysis is a non-symmetric diffusion
determined by $A(\msx)=A_0(\msx)+B(\msx)$, where $B(\msx)=\sum_i b_i(\msx)
\partial_i$, and a construction of its approximations by MJPs in lattices. For
$d=2$ the method is valid for any $a$, and for $d>2$ the method is valid
for tensor valued functions $a$ satisfying an additional constraint.
In this way we offer an efficient method
for simulation of sample paths of a non-symmetric generalized diffusion
using MJPs. Each $X_n(\cdot)$ can be simulated by well-known methods.

In the case of classical diffusion determined by an elliptic operator of the
form $A(\msx) = -\sum_{ij} a_{ij}(\msx) \partial_i
\partial_j + \sum_i b_i(\msx) \partial_i$, where $a_{ij}, b_i$ are H\"{o}lder
continuous on ${\bbR}^d$, approximations by MJPs
can be efficiently used to simulate the first exit from a bounded set
of ${\bbR}^d$. However, such an approach is one of several
existing possibilities, and the motivation for the construction in terms of MJPs
happens to be of lesser importance, since the process has a representation in
terms of SDE which can be simulated straightforwardly. For one-dimensional
generalized diffusion, there exist representations in terms of SDE as described by
\'Etor\'e \cite{Et} and Lejay \& Martinez \cite{LM}, so that such
representations can be used for simulation. In the case
of a process defined by a differential operator in divergence form on
${\bbR}^d$, $d\geq 2$, there is no such natural representation, so approximations
by MJPs are essentially the only tool available for simulations.

A class of convergent approximations $X_n(\cdot)$ in~\cite{SZ} is constructed
by using discretizations of the corresponding Dirichlet form $(v,u) \mapsto
a(v,u)$. The functions $\partial_iv,\partial_iu$ in $a(v,u)$ are
approximated by forward difference operators in local basis which generally
varies. This approach is anticipated in~\cite{MW} without any remarks on the
convergence of the constructed MJPs $X_n(\cdot)$. In our approach, the MJPs
$X_n(\cdot)$ are constructed in terms of generators $A_n({\rm gen})$ on the grids,
\begin{equation}\label{exp1.1}
 G_n \ = \ \Big\{h_n\,\sum_{i=1}^d \,k_i\mse_i\,:\, k_i \in {\bbZ}\Big\},\quad
 h_n = 2^{-n},
\end{equation}
with a fixed basis $\{\mse_i\}_1^d$. The index set of grid-knots is denoted by
$I_n$. In order to simplify expressions, we often write $h$ instead of $h_n$.
The generators $A_n({\rm gen})$ are constructed explicitly from a general
principle which is not directly related to forward difference operators.
Then the discretizations $h^da_n(v,u)$ of the original Dirichlet form are
associated to the constructed generators $A_n({\rm gen})$. It turns out that
$h^da_n(v,u)$ cannot be simply obtained from discretizations of original form
by using forward/backward difference operators. The obtained class of
$A_n({\rm gen})$ is not included among generators constructed in~\cite{SZ}.

Advantages and drawbacks of the present construction can be briefly
described as follows.

{\em Advantages:} The convergence of MJPs is proved for a non-symmetric
generalized diffusion. The generators $A_n({\rm gen})$ are explicitly
given in terms of values of functions $a_{ij}, b_i$.
For the case of $d=2$, the matrix entries of $A_n({\rm gen})$
can be easily implemented into a computer code because rotations of
coordinates are avoided. For $d=2$, the construction
holds for a general matrix valued function $a(\cdot)$ on ${\bbR}^2$.

{\em Disadvantages:} For $d \geq 3$, the proposed construction is not valid
for all $a(\cdot)$ on ${\bbR}^d$. The restriction is defined by an auxiliary
matrix valued function $\hat{a}(\cdot)$, with the following matrix entries:
\begin{equation}\label{exp1.2}
 \hat{a}_{ii} \ = \ a_{ii}, \quad \hat{a}_{ij} \ = \ -|a_{ij}|, \ i \ne j.
\end{equation}
The here proposed construction of $A_n({\rm gen})$ is valid only if $\hat{a}(\cdot)$
is strictly positive definite on ${\bbR}^d$. For $d \geq 3$, there are simple
examples of pairs $a, \hat{a}$, where $a$ is definite and $\hat{a}$ is indefinite.
We can say that (\ref{exp1.2})
is valid if the off-diagonal entries of $a(\cdot)$ are small in comparison with
the diagonal ones. When $\hat{a}(\msx)$ becomes indefinite, it is necessary to
apply a local rotation of coordinates, that ensures a diminishing of the
off-diagonal entries, thus ensuring the positive definiteness of $\hat{a}(\msx)$
in the new coordinate system.

Now we can describe basic steps in the construction and the proofs.
Let $U(\cdot)$ be the strongly continuous semigroup in the Banach space
$\dot{C}({\bbR}^d)$ (continuous functions vanishing at infinity)
which is associated with the diffusion process $X(\cdot)$ \cite{EK}.
For each $n$ the space of discretizations of $\dot{C}({\bbR}^d)$
in terms of grid-functions on $G_n$ is denoted by $\dot{l}_\infty(G_n)$.
The generators $A_n({\rm gen})$ determine semigroups $U_n(\cdot)$ in
$\dot{l}_\infty(G_n)$. There exist continuous mappings
$\Phi_n: \dot{l}_\infty(G_n) \mapsto \dot{C}({\bbR}^d)$
and $\Phi_n^{-1}:\dot{C}({\bbR}^d) \mapsto \dot{l}_\infty(G_n)$ with the
following properties: $\Vert \Phi_n \Vert =1, \Vert \Phi_n^{-1} \Vert \leq 1$
and $\Phi_n^{-1}\Phi_n=I$ on $\dot{l}_\infty(G_n)$. For $f \in \dot{C}({\bbR}^d)$,
the grid-function ${\bf f}_n=\Phi_n^{-1}f$ is called the discretization of
$f$. Now, for each $t \geq 0$, the grid-functions
${\bf u}_n(t)=\Phi_n^{-1} U(t)f$ and $U_n(t){\bf f}_n$ can be compared. We
need the following relation:
\begin{equation}\label{exp1.5}
 \lim_{n \to \infty}\:\sup\{ \lnorm U_n(t){\bf f}_n - \Phi_n^{-1}
 U(t)f\lnorm_\infty\,:\,t\in [0,1]\} \ = \ 0,
\end{equation}
where $\lnorm \cdot\lnorm_\infty$ is the $l_\infty(G_n)$-norm. By Theorem 2.11
of \cite{EK}, this relation is a sufficient condition for the convergence
in distribution of MJPs $X_n(\cdot)$ to the diffusion process $X(\cdot)$.

Definitions of various objects are given in Section \ref{sec2}. In Section
\ref{sec4}, a class of generators $A_n({\rm gen})$ is constructed explicitly.
Grid-functions ${\bf u}_n(t)=\exp(A_n({\rm gen})t){\bf f}_n$ and their images
$u(n,t)=\Phi_n {\bf u}_n(t)$ in $\dot{C}({\bbR}^d)$ are studied in
Section \ref{sec5}. The convergence $u(n,t)\mapsto U(t)f$ in $\dot{C}({\bbR}^d)$
is proved initially for the case of smooth coefficients $a_{ij}, b_i$, and then for
general coefficients $a_{ij}, b_i$. A numerical example is given in the last
section.

\section{PRELIMINARIES}\label{sec2}
The Euclidean norm in ${\bbR}^d$ is denoted by $|\cdot|$.
All the open subsets of ${\bbR}^d$, considered in this work,
are bounded and connected open sets with Lipshitz boundary \cite{Ma,Se}.
We call a subset of this kind a {\em Lipshitz domain}, denote it by D, and
its boundary by $\partial D$.

The Banach spaces of functions $C^{(k)}({\bbR}^d), C_0^{(k)}({\bbR}^d)=
C_0({\bbR}^d)\cap C^{(k)}({\bbR}^d)$ are defined as usual, $C_0({\bbR}^d)$ being
the linear space of continuous functions with compact support.
Their norms are denoted by $\Vert \cdot \Vert_\infty^{(k)}$.
The closure of functions in $C({\bbR}^d)$ with compact support determines the
Banach space $\dot{C}({\bbR}^d)$. The H\"{o}lder space of parameter
$k+\alpha, k\in {\bbN}_0, \alpha \in (0,1)$, is denoted by $C^{(k+\alpha)}({\bbR}^d)$
and defined as the completion of $C_0^{(\infty)}({\bbR}^d)$ in the norm:
\[\Vert u \Vert^{(k+\alpha)} \ = \ \Vert u\Vert_\infty^{(k)}\,+\,\sup \Big\{
 \frac{|\partial^k u(\msx+\msz)-\partial^k u(\msx)|}{|\msz|^\alpha} \::\:
 \msx \in {\bbR}^d,  0<|\msz|\leq 1 \Big\}.\]
The $L_p$-spaces as well as Sobolev $W_p^1$-spaces are defined in a
standard way \cite{Ma,Se}. Their norms are denoted by
$\Vert \cdot \Vert_p$ and $\Vert \cdot \Vert_{p,1}$, respectively. For each
$p,\: 1 \leq p \leq \infty$, the norm of $W_p^1({\bbR}^d)$ is defined
by $\Vert u \Vert_{p,1} =  \big(\Vert u \Vert_p^2 + \Vert \nabla u
\Vert_p^2\big)^{1/2}$, where $\Vert \,\nabla u\, \Vert_p \, = \, \big(\sum_{j=1}^d\,
\Vert \,\partial_j u\, \Vert_p^2\big)^{1/2}$. We say that a function $f$ on
${\bbR}^d$ belongs to a class $C^{(k+\alpha)}$ on ${\bbR}^d$ if
$\Vert f\Vert^{(k+\alpha)}<\infty$.

In this article we consider a $2^{{\rm nd}}$-order elliptic operator on ${\bbR}^d$,
\begin{equation}\label{exp3.1}
 A(\msx) = -\:\sum_{ i,j=1}^d  \partial_i a_{ij}(\msx) \partial_j
 \:+\: \sum_{j=1}^d b_j(\msx)\, \partial_j,
\end{equation}
for which the coefficients fulfill the following:

\begin{assumption}\label{Ass3.1}
The functions  $a_{ij}= a_{ji}$, $b_{i}, i,j = 1,2,\dots ,d$,
are measurable on ${\bbR}^d$ and have the following properties:
\begin{description}\itemsep 0.cm
 \item{a)} There exist positive numbers $\Md, \Mg,\:0 < \Md \leq \Mg$,
such that the strict ellipticity is valid:
\begin{equation}\label{exp3.2}
 \Md \,|\msz|^2  \,\leq\,\sum_{i,j=1}^{d}  a_{ij}(\msx)z_i \bar{z}_{j}
 \leq \,\Mg\, |\msz|^2 ,  \quad  \msx, \msz \in  {\bbR}^d.
\end{equation}
 \item{b)} The functions $b_i$ are bounded on ${\bbR}^d$.
\end{description}
\end{assumption}

In our analysis we regularly use the notation $A_0(\msx)=-\sum_{ ij=1}^d
\partial_i a_{ij}(\msx) \partial_j, B(\msx) =\sum_{j=1}^d b_j(\msx) \partial_j$ and
$A(\msx)=A_0(\msx)+B(\msx)$. The following real bilinear form on
$W_q^1(D) \times W_p^1(D)$, $1/p+1/q = 1$, defined by:
\begin{equation}\label{exp3.4}
 a(v,u) \ = \ \sum_{i,j =1}^d\: \int\: a_{ij}( \msx)\:
 \partial_i v(\msx)\:  \partial_j u(\msx) \:d\msx
 - \sum_{i=1}^d \: \int\: b_i(\msx) \:v(\msx)\,
 \partial_i u(\msx) \:d \msx
\end{equation}
is associated with the differential operator $A(\msx)$.

A basic result towards a proof of (\ref{exp1.5}) is the following theorem \cite{St}:

\begin{theorem}\label{Th3.2} Let the differential operator $A(\msx)$ be defined
by (\ref{exp3.1}) and Assumption \ref{Ass3.1}. Then $-A(\msx)$ has the closure
in $\dot{C}({\bbR}^d)$ and generates a Feller semigroup in $\dot{C}({\bbR}^d)$.
\end{theorem}

The semigroup of this theorem has the restriction to a closed subspace of
$\dot{C}^{(\alpha)}({\bbR}^d)=\dot{C}({\bbR}^d)\cap C^{(\alpha)}({\bbR}^d)$
which is also a strongly continuous semigroup.

\begin{corollary}\label{Cor3.1} Let $A(\msx)$ be defined as in Theorem
\ref{Th3.2}. There exists an $\alpha \in (0,1)$ and $\sigma\geq 0$ that depend only
on $\Md, \Mg$ and $\Vert \msb \Vert_\infty$, such that the following two
assertions are valid:
\begin{description}\itemsep -0.1cm
 \item{(i)} The operators $U(t)$ in $\dot{C}^{(\alpha)}({\bbR}^d)$ are bounded
uniformly on segments of $[0,\infty)$.
 \item{(ii)} There exists a closed subspace $F^{(\alpha)}({\bbR}^d)\subseteq
\dot{C}^{(\alpha)}({\bbR}^d)$ such that the closure of $-A(\msx)$ in $F^{(\alpha)}
({\bbR}^d)$ generates a strongly continuous semigroup $U(\cdot)$ in
$F^{(\alpha)}({\bbR}^d)$, with $\Vert U(t)\Vert_\infty^{(\alpha)}\leq
\exp(\sigma t)$.
\end{description}
\end{corollary}

{\Proof} The Feller semigroup of Theorem \ref{Th3.2} can be represented as
$t \mapsto \exp(t)V(t)$, where $V(\cdot)$ is the strongly continuous semigroup
generated by the closure of $-(I+A(\msx))$ in $\dot{C}({\bbR}^d)$.
Due to Theorem \ref{Th3.2} the differential
operator $I+A(\msx)$ has the bounded inverse $(I+A)^{-1}$ in $\dot{C}({\bbR}^d)$.
For any $f\in \dot{C}({\bbR}^d)$ the continuous function $(t,\msx)\mapsto
u(t,\msx)$, defined by $u(t)=V(t)f$, is a solution to the following IVP for PDE
\begin{equation}\label{exp3.5}
 \partial_t u(t,\msx) +(I +A(\msx))u(t,\msx) \ = \ 0,
\end{equation}
with the initial condition $u(0,\msx)=f(\msx)$. In addition,
$\Vert u(t)\Vert_\infty \leq \Vert f\Vert_\infty$ for all $t\geq 0$.
It turns out that a solution to (\ref{exp3.5}), for each $t>0$, is an element of
$\dot{C}^{(\alpha)}({\bbR}^d)$. For an initial condition $f\in \dot{C}^{(\alpha)}
({\bbR}^d)$ we have a stronger result, a solution to (\ref{exp3.5}) is
an element of $\dot{C}^{(\alpha)}({\bbR}^d)$ uniformly with respect to $t\in
[0,1]$. This follows from the following arguments. Let us consider
the balls $B_1(\msx)\subset {\bbR}^d$, centered at $\msx\in {\bbR}^d$ and having
the radius equal to 1, and let us consider the restrictions $u(t)|B_1(\msx)$
for any $\msx\in {\bbR}^d$ and any $t\in [0,1]$. According to the basic
theorem of Section 10, Chapter 3 of \cite{LSU}, there
exists $\alpha\in (0,1)$ depending on $\Md, \Mg, \Vert \msb \Vert_\infty$ and
$\beta$ depending on $\Md, \Mg, \Vert \msb \Vert_\infty$ and $\Vert f
\Vert_\infty^{(\alpha)}$ such that a solution
$u$ to (\ref{exp3.5}) on $B_1(\msx)$ is an element of $C^{(\alpha)}(B_1(\msx))$, and
\begin{equation}\label{exp3.3}
 \sup_{t\in [0,1]}\,\Vert u(t)|B_1(\msx) \Vert^{(\alpha)} \,+\,
 \sup_{t,s\in [0,1]}\,\sup_{\msyd\in B_1(\msxd)}\,
 \frac{| u(t,\msy)-u(s,\msy)|}{|t-s|^\alpha} \ \leq \ \beta.
\end{equation}
This inequality implies that the operators $V(t), t\in [0,1]$, are linear in
$\dot{C}^{(\alpha)}({\bbR}^d)$, with a norm which is uniformly bounded with
respect to $t\in [0,1]$. Hence, the operators $V(t)$ define a semigroup of
bounded operators in $\dot{C}^{(\alpha)}({\bbR}^d)$. This proves the assertion
(i) of the corollary.

For any $g\in \dot{C}({\bbR}^d)$ the function $f=(I+A)^{-1}g$ belongs to the
domain $\mathfrak{D}(A)$ of the closure of $A(\msx)$ in $\dot{C}({\bbR}^d)$,
and $V(t)f=(I+A)^{-1}V(t)g$. The linear space spanned by $f=(I+A)^{-1}g,\: g \in
\dot{C}({\bbR}^d)$, is denoted by $E({\bbR}^d)$.
Let us assume now that we can prove $f=(I+A)^{-1}g \in \dot{C}^{(\alpha)}
({\bbR}^d)$ for any $g \in \dot{C}({\bbR}^d)$, i.e. $\Vert f\Vert^{(\alpha)}
\leq \beta \Vert g\Vert_\infty$, where $\beta$ is a number depending on $\Md,
\Mg$ and $\Vert \msb \Vert_\infty$. This would imply that $E({\bbR}^d)$
consists of elements in $\dot{C}^{(\alpha)}({\bbR}^d)$ and its completion in the
$\Vert \cdot \Vert^{(\alpha)}$-norm is a closed space $F^{(\alpha)}({\bbR}^d)
\subseteq \dot{C}^{(\alpha)}({\bbR}^d)$.
In addition, the above assumption would imply the following inequality:
\begin{equation}\label{exp3.5b}
 \Vert (V(t)-I)f \Vert^{(\alpha)} \ \leq \ \beta \,\Vert  (V(t)-I) g \Vert_\infty,
\end{equation}
i.e. the continuity of the function $t\mapsto V(t)$ on a dense linear space
$E({\bbR}^d)\subset F^{(\alpha)}({\bbR}^d)$. Since the operators
$V(t)$ are bounded in $\dot{C}^{(\alpha)}({\bbR}^d)$ uniformly with respect to $t\in
[0,1]$, the obtained inequality would imply that the operators $V(\cdot)$ define
a strongly continuous semigroup in the Banach space $F^{(\alpha)}({\bbR}^d)$.
Then the semigroup $U(\cdot)$ of (ii) is also strongly continuous on $F^{(\alpha)}
({\bbR}^d)$ and the existence of $\sigma\geq 0$ in (ii) follows from the theory
of strongly continuous semigroups.

Therefore, in order to prove (ii) it remains to show
that $(I+A)^{-1}$ maps $\dot{C}({\bbR}^d)$ into
$\dot{C}^{(\alpha)}({\bbR}^d)$. For this we consider the elliptic problem:
\begin{equation}\label{exp3.5c}
 \big(I\,+\,A(\msx)\big)\,u(\msx) \ = \ f(\msx), \quad \msx \in {\bbR}^d,
\end{equation}
where $f\in \dot{C}^{(\alpha)}({\bbR}^d)$. The existence of solution $u=(I+A)^{-1}f$
in $\dot{C}({\bbR}^d)$ is granted by Theorem \ref{Th3.2}, and in addition we have
$\Vert u\Vert_\infty \leq \Vert f\Vert_\infty$. In order to prove that
$u\in \dot{C}^{(\alpha)}({\bbR}^d)$, it is sufficient to verify that $u|B_1(\msx)\in
C^{(\alpha)}(B_1(\msx))$ for any $\msx$ and a fixed $\alpha$. To obtain this,
we apply the theorem in Section 14, Chapter 3 of \cite{LU}. The following estimate
is valid:
\[ \Vert u| B_1(\msx) \Vert^{(\alpha)} \ \leq \ \beta\,\Vert f\Vert_\infty, \]
where both, $\alpha\in (0,1)$ and $\beta$, depend on $\Md, \Mg, \Vert \msb
\Vert_\infty$.

The parameters $\alpha \in (0,1)$ for the parabolic (\ref{exp3.5}) and elliptic
problems (\ref{exp3.5c}) are not necessarily equal. We choose the minimum
to obtain (\ref{exp3.5b}), and consequently the assertion (ii). {\QED}

The assertion (i) of Corollary \ref{Cor3.1} is used in the proof of Theorem
\ref{Th5.2}, which establishes the basic inequality (\ref{exp1.5}).

The linear space of grid-functions on $G_n$ is denoted by $l(G_n)$.
Elements of $l(G_n)$ are also called {\em columns}. Columns are denoted by
${\bf u}, {\bf v}$ etc, while their entries are denoted by $u_{\mskd}, v_{\msld}$
etc. Thus a column ${\bf u}_n$ has entries $({\bf u}_n)_{\mskd}, \msk\in I_n$.
Columns with a finite number of components span a linear space $l_0(G_n)$.
The completion of $l_0(G_n)$ in the $\lnorm \cdot \lnorm_\infty$-norm is denoted
by $\dot{l}(G_n)$. The corresponding $l_p$-spaces
are denoted by $l_p(G_n)$ and their norms by $\lnorm \cdot \lnorm_p$.
For $p<\infty$ this norm is defined by $\lnorm {\bf u} \lnorm_{p} =
\big[\sum_{\mskd \in I_n}\: |u_{\mskd}|^p \big]^{1/p}$,
and for $p = \infty$ by $\lnorm {\bf u} \lnorm_{\infty} = \sup\{|u_{\mskd}|: \msk
\in {\bbZ}^d\}$. The scalar product in
$l_2(G_n)$ is denoted by $\lev \cdot |\cdot \des$ and sometimes by $(
\cdot |\cdot )$.

Let $f \in C({\bbR}^d)$ and the column ${\bf f}_n \in l(G_n)$ be defined by
its components $({\bf f}_n)_{\mskd}=f(h\msk)$, where $h=h_n=2^{-n}$. Then the
mapping $C({\bbR}^d) \owns f \mapsto {\bf f}_n \in l(G_n)$ is called the
discretization of $f$.

The shift operator $Z(\msx), \msx \in {\bbR}^d$, acting on functions
$f:{\bbR}^d \mapsto {\bbR}$, is defined by $\big(Z(\msz)f\big)(\msx) =
f(\msx+\msz)$. Similarly we define the discretized shift operator. The shift
operator $Z_n(r,i)$ by $r$ units in the direction $\mse_i$ is defined by
$\big ( Z_n(r,i) {\bf u}_n \big)_{\mskd} = ({\bf u}_n)_{\msld}$, where $\msl =
\msk + r\mse_i$. The partial derivatives of $u \in C^{(1)}({\bbR}^d)$
with respect to a grid step $h$ are discretized by
forward/backward finite difference operators in the usual way,
\begin{equation}\label{ex2.7}\begin{array}{c}
 \der_i(th) u(\msx) \ = \ \frac{1}{th} \big(
 u(\msx + th \mse_i) \:-\: u(\msx) \big ),\\
 \widehat{\der}_i(th) u(\msx) \ = \ \frac{1}{th} \big(
 u(\msx) \:-\: u(\msx - th\mse_i) \big ), \end{array} \quad \msx \in {\bbR}^d, \ t
 > 0.
\end{equation}
Let $r \in {\bbZ}\setminus \{0\}$.
Discretizations of the functions $\partial_iu$ on $G_n$, denoted by
$U_i(r){\bf u}_n, V_i(r){\bf u}_n$, are defined by:
\[ \big (U_i(r)\,{\bf u}_n\big )_{\msmd} \:=\: \der_i\big(rh \big)
 \:u(\msx_{\msmd}), \quad  \big (V_i(r)\,{\bf u}_n\big )_{\msmd} \:=\:
 \widehat{\der}_i \big(rh \big) \:u(\msx_{\msmd}),\]
where $\msx_m \in G_n$. Then
\begin{equation}\label{ex2.8}
 \begin{array}{l} U_i(r) \ = \ (rh)^{-1}(Z_n(r,i) \,-\, I),\\
 V_i(r) \ = \  (rh)^{-1} \big(I - Z_n(-r,i) \big) \:= \:
 U_i(-r) \:=\:\,-\, U_i(r)^T.
\end{array}
\end{equation}
Therefore we have $U_i(-r) = U_i(r)\,Z_n(-r,i) = Z_n(-r,i)\,U_i(r)$, and
similarly for $V_i(\cdot)$. We use the abbreviations $U_i=U_i(1), V_i=V_i(1)$.

A matrix $A_n$ on $G_n$ is called a matrix of {\em positive type} if the diagonal
entries of $A_n$ are positive, off-diagonal entries are non-positive and the row
sums are non-negative. If $A_n(gen)$ is the generator of a MJP in $G_n$, then
$-A_n(gen)$ is a matrix of positive type.

\section{CONSTRUCTION OF GENERATORS OF MJPs}\label{sec4}
To discretize $A(\msx)$ means to associate to $A(\msx)$ a sequence of matrices
$A_n$ on $G_n, n \in {\bbN}$, with the following properties:
\[ a(v,u) \ = \ \lim_{n\to\infty}\:h^d\lev {\bf v}_n|A_n{\bf u}_n\des,
 \quad v,u \in C_0^{(1)}({\bbR}^d).\]
The terminology "discretizations" of $A(\msx)$ instead of approximations of
$A(\msx)$ appears to be more suitable at the beginning of the construction, since
the convergence analysis is postponed to Section \ref{sec5}.

We wish to emphasize that discretizations $A_n$ of the differential operator
$A_0(\msx)$ are derived from a general principle, similar to the one exploited in
\cite{LR2}. This method is not based on finite difference formulas. Nevertheless,
bilinear forms need to be associated to $A_n$ so that $A_n$ can be derived from
the corresponding variational equalities. The constructed bilinear forms can be
considered as the discretizations of the original form (\ref{exp3.4}). These
forms are basic objects in our proof of convergence of discretizations.

Discretizations to be considered in this section are possible if certain conditions
on $a_{ij}$ are fulfilled. The required conditions are stronger than those given in
Assumption \ref{Ass3.1}. By relaxing them gradually as $n \to \infty$ we obtain
discretizations for a general $A(\msx)$ given in Assumption \ref{Ass3.1}.

To each pair $\msv \in G_n, \msp \in {\bbN}^d$, there is associated a
rectangle $C_n(\msp,\msv) =  \prod_{i=1}^d  [v_i,v_i +hp_i)$ with the "lower left"
vertex $\msv$ and the edge of size $h p_i$ in the $i$-th coordinate direction.
These rectangles define a partition of ${\bbR}^d$.
Apart from these rectangles, we will need the closed rectangles,
\begin{equation}\label{exp4.s}
 S_n(\msp,\msv) \ = \ \prod_{i=1}^d \: [v_i - hp_i\,,\,v_i+hp_i],
\end{equation}
which are centered at the grid-knots $\msv$. Evidently, $S_n(\msp,\msv)$
is the union of closures of those rectangles $C_n(\msp,\msx)$
which contain the grid-knot $\msv$.

A discretization $A_n$ is defined in terms of its matrix entries $(A_n)_{\mskd\msld}$,
where $h\msk, h\msl \in G_n$. For a fixed $\msx = h\msk \in G_n$ the set of all
the grid-knots $\msy = h\msl$ such that $(A_n)_{\mskd\msld} \ne 0$ is denoted by
${\cal N}(\msx)$ and called the {\em numerical neighborhood} of $A_n$ at $\msx \in
G_n$. The set ${\cal N}(\msx)$ contains always a subset consisting of $\msx$
and $2d$ elements $\msx\pm h\mse_i, i=1,2,\ldots,d$. Additional elements of
${\cal N}(\msx)$ are possible as we shall see, depending on the sign of
$a_{ij}, i \ne j$. In terms of the MJP $X_n(\cdot)$, the set ${\cal N}(\msx)$
consists of the states of possible jumps from the state $\msx$.
Let us point out that the sets ${\cal N}(\msx)$ vary with $n$, that is
for two grids $G_n, G_m, n\ne m$ and $\msx\in G_n\cap G_m$, the
corresponding numerical neighborhoods ${\cal N}(\msx)$ are different.

\subsection*{General setup}
In order to give a comprehensive insight into the proposed construction of
the discretizations, it is convenient to initially consider a differential
operator $A_0=-\sum_{ij}a_{ij}\partial_i \partial_j$ having a constant diffusion
tensor $a=\{a_{ij}\}_{11}^{dd}$. In this case, the matrices $A_n$ need to have a
property which is called the {\em consistency}. Let $\msx\mapsto p(\msx)=p_0+\sum_i
p_ix_i+\sum_{ij} p_{ij}x_i x_j$, $p_0,p_i, p_{ij}\in \bbR$, be a second degree
polynomial in arguments $x_i$, and let ${\bf p}_n$ be its
discretization on grid $G_n$. Then the consistency holds if the following
identities are valid:
\begin{equation}\label{exp4.0}
 \delta(A,n,\msx)p(\msx) \ := \ A(\msx)p(\msx)-\big(A_n {\bf p}_n\big)_{\mskd}
 \ = \ 0,\quad \msx=h\msk \in G_n.
\end{equation}
These consistency conditions are sufficient for proving the convergence of
$U_n(t){\bf f}_n$ to $\Phi_n U(t)f$ in the Banach space of continuous functions on
${\bbR}^d$ as specified in (\ref{exp1.5}). Actually,
the first step of construction of matrices $A_n$ begins with a search for
those matrices $A_n$ which are simultaneously of positive type and fulfill
the conditions $\delta(A,n,\msx)p(\msx)=0$ on $G_n$ \cite{LR1}.

We request that the discretizations $A_n$ have the following properties:
\begin{description}\itemsep -0.1cm
 \item{a)} The numerical neighborhoods ${\cal N}(\msx)\subset G_n$ resemble
each other, that is ${\cal N}(\msx)=\msx+{\cal N}({\bf 0})$.
 \item{b)} The numerical neighborhoods ${\cal N}({\bf 0})\subset G_n, n \in {\bbN}$,
resemble each other, i.e. if $\msy =h_n\msl \in {\cal N}({\bf 0})\subset
G_n$ then $\msy' =h_m\msl \in {\cal N}({\bf 0})\subset G_m$, $n,m \in \bbN, n\ne m$.
 \item{c)} The matrices $A_n$ are symmetric.
 \item{d)} The matrices $A_n$ are consistent discretizations of $A_0$.
\end{description}
Now we have the following result.

\begin{lemma}\label{lem2.1} If $\hat{a}$ is positive definite then there exist
matrices $A_n$ of positive type fulfilling the conditions a)-d). The non-trivial
entries of $A_n$ are defined in terms of $d$ natural numbers $r_1, r_2,\ldots,
r_d\in \bbN$ by the following expressions:
\begin{equation}\label{exp2.13}\begin{array}{llll}
 \big(A_n\big)_{\mskd \,\mskd} &=& -\,\sum_{h\msld \in {\cal N}(h\mskd), \
 \msld \ne \mskd} \big(A_n\big)_{\mskd \msld},&\\
 \big(A_n\big)_{\mskd \,\mskd\pm \msed_i} &=& -\,\frac{1}{h^2}\,
 \left [a_{ii}-\sum_{m\ne i}\frac{r_i}{r_m}|a_{im}|\right ],& i=1,2,\ldots,d,\\
 \big(A_n\big)_{\mskd \,\mskd\pm \mszd(i,j,\pm)} &=&-\, \frac{1}{h^2 r_ir_j}\,
 |a_{ij}|,& i,j=1,2,\ldots,d,\end{array}
\end{equation}
where in the last line, $\msz(i,j,+) = r_i\mse_i + r_j\mse_j$ corresponds to the
case $a_{ij}>0$, and $\msz(i,j,-)=r_i\mse_i-r_j\mse_j$ corresponds to the case
$a_{ij}<0$.
\end{lemma}

{\Proof} For a matrix $A_n$ defined by (\ref{exp2.13}) the conditions a)- c)
are obviously satisfied. (Recall that ${\cal N}(h\msk)$ consists of those
grid-knots $h\msl \in {\cal N}(h\msk)$ for which $\big(A_n\big)_{\mskd \,\msld}$
are non-trivial). It remains to prove the condition d), and that $A_n$ is
of positive type if $\hat{a}$ is a positive definite matrix.

The condition d) is equivalent to the following property.
Let $\msx \mapsto p(\msx)=\sum_{ij}p_{ij}x_ix_j +\sum_i p_ix_i+p_0$ be a second
degree polynomial. Then d) is valid if and only if
\[\sum_{h\mskd \in {\cal N}({\bf 0})}\:\big(A_n\big)_{{\bf 0}\mskd}\,p(h\msk) \ = \ -2\,
 \sum_{i \leq j} a_{ij}p_{ij}.\]
For $A_n$ defined by (\ref{exp2.13}) this identity can be easily verified
by using the monomials $p(\msx)=x_ix_j$ with $i \ne j$ and $p(\msx)=x_i^2$, for
all $i,j = 1,2,\ldots,d$.

In the last step of the proof we show that there exists a sequence
$r_1, r_2,\ldots,r_d\in \bbN$ for which the brackets in the second line of
(\ref{exp2.13}) are positive, thus ensuring the matrix $A_n$ to be of positive
type. For this purpose we assume that the matrix $\hat{a}$ is positive definite.
Let us consider the eigenvalue problem $\hat{a}\msw =\lambda \msw$, and let
us also assume that the matrix $\hat{a}$ is irreducible.
For $\mu >0$ sufficiently large the irreducible matrix
$\mu I+\hat{a}$ has the inverse $(\mu I+\hat{a})^{-1}$ with positive entries,
so the Perron-Frobenius theorem can be applied to $(\mu I+\hat{a})^{-1}$.
Thus the eigenvector corresponding to its maximal eigenvalue is positive. This
result can be also formulated in terms of the problem $\hat{a}\msw =\lambda_1 \msw$
for the minimal eigenvalue $\lambda_1$ of $\hat{a}$. We have $\lambda_1 \msw>0$ and
consequently $\hat{a}\msw >0$. This inequality can be rewritten as
$a_{ii}-\sum_{m\ne i}(q_i/q_m)|a_{im}|>0$ where $q_i=w_i^{-1}$. One can find
rational approximations
$r_i/r_m$ of $q_i/q_m$ which also fulfil the obtained inequalities.
If $\hat{a}$ is
not irreducible, it can be rewritten in a block diagonal matrix form with irreducible
blocks. The previous construction can be applied to each block. {\QED}

In the case of $d=2$, either both $a,\hat{a}$ are positive definite or neither
is. For $d=3$, there are symmetric positive definite matrices $a$
for which $\hat{a}$ are indefinite. For instance, the symmetric matrix $a$
of order $d=3$ defined by $a_{ii} = 1, a_{12} = a_{23} = -1/\sqrt{2}$ has
positive eigenvalues for the case of $a_{13} > 0$ and a negative eigenvalue for
the case of $a_{13} < 0$.

For $d=2$ and $r_1 = 3, r_2 =1$, two possible
numerical neighborhoods ${\cal\ N}(\msx)$ are illustrated in
Figure~\ref{fig4.2}.

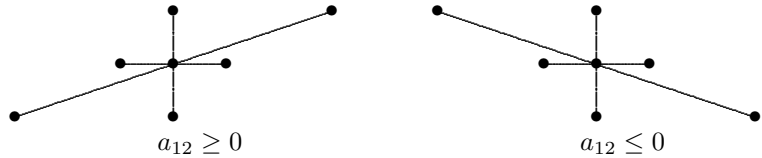
\begin{Figure}\label{fig4.2}
\begin{center}
\beginpicture
\setcoordinatesystem units <1.0pt,1.0pt>
\setplotarea x from -40 to 300, y from 20 to 100
\setlinear
\plot 40 60, 80 60 /
\plot 60 40, 60 80 /
\plot 0 40, 120 80 /
\plot 200 60, 240 60 /
\plot 220 40, 220 80 /
\plot 160 80, 280 40 /
\put{$\bullet$} [c] at 60 60
\put{$\bullet$} [c] at 40 60
\put{$\bullet$} [c] at 80 60
\put{$\bullet$} [c] at 60 40
\put{$\bullet$} [c] at 60 80
\put{$\bullet$} [c] at 0 40
\put{$\bullet$} [c] at 120 80
\put{$\bullet$} [c] at 220 60
\put{$\bullet$} [c] at 200 60
\put{$\bullet$} [c] at 240 60
\put{$\bullet$} [c] at 220 40
\put{$\bullet$} [c] at 220 80
\put{$\bullet$} [c] at 160 80
\put{$\bullet$} [c] at 280 40
\put{$r_1 = 3, \ r_2 =1$} [c] at 150 5
\put{$a_{12} \geq 0$} [c] at 70 30
\put{$a_{12} \leq 0$} [c] at 230 30
\endpicture
\end{center}
\centerline{Numerical neighborhoods}
\end{Figure}

If $A_n$ is a matrix from Lemma \ref{lem2.1}, then $A_n(gen)=-A_n$ is a generator
of a MJP in $G_n$. Let us consider now functions $a_{ij}$ on ${\bbR}^d$ defining
a diffusion tensor at each $\msx \in {\bbR}^d$. Let us assume in addition that
the matrix $\hat{a}(\msx)=\{\hat{a}_{ij}(\msx)\}_{11}^{dd}$ is positive definite
at each $\msx \in {\bbR}^d$. One could replace the numbers $a_{ij}$ of (\ref{exp2.13})
by the numbers $a_{ij}(h\msk)$. The resulting $A_n(gen)=-A_n$ would again be a
generator of a MJP in $G_n$. However, the matrices $A_n$ thus obtained
are discretizations of the
differential operator $A(\msx)=-\sum_{ij}a_{ij}(\msx)\partial_i \partial_j$, as will
be seen in Section \ref{sec5}, and not of $A(\msx)=-\sum_{ij}\partial_i a_{ij}(\msx)
\partial_j$. Here we present a method of construction of discretizations of
a differential operator in divergence form, $A(\msx)=-\sum_{ij}\partial_i a_{ij}
(\msx)\partial_j$, resulting in matrices of positive type, with a structure
similar to (\ref{exp2.13}). The operator in divergence form, $A(\msx)$, naturally
corresponds to the bilinear form $a(v,u)=\sum_{ij}\int a_{ij}(\msx)\partial_iv(\msx)
\partial_ju(\msx)d\msx$. Consequently, its discretizations $A_n$ will correspond to
a sequence of discretized forms $a_n(v,u)$ on grids $G_n$.

\subsection*{Construction for $d=2$}
Initially one considers the
constant coefficients $a_{ij}$ and constructs the forms $a_n(v,u)=
\lev {\bf v}_n|A_n{\bf u}_n\des$ for which the matrices $A_n$ coincide with
(\ref{exp2.13}). Then the obtained expression of $a_n(v,u)$ is generalized to
the case of non-constant coefficients $a_{ij}$. The corresponding matrices $A_n$
are obtained from the variational method in the standard way.

For the case of constant coefficients $a_{ij}$, the bilinear form $h^2a_n(v,u)$ that
discretizes the form $a(v,u)=\sum_{i,j=1}^2a_{ij}\partial_iv \partial_ju$, will be
a second order polynomial in the quantities $\der_i(r_ih)v(\msx)$, $\der_j(r_jh)u
(\msx)$ with a certain choice of $r_i, r_j\in \bbN$ and $\msx \in G_n$. In order to
write down the form as simply as possible, we use the following abbreviations:
\[\begin{array}{lllll}
 u_{\mskd}(i,r_i)&=&\der_i(r_ih)u(h\msk)&=& (r_ih)^{-1}\,[u(h\msk+r_ih\mse_i)
 -u(h\msk)],\\
 \widehat{u}_{\mskd}(i,r_i)&=&\widehat{\der}_i(r_ih)u(h\msk)&=& (r_ih)^{-1}\,
 [u(h\msk) -u(h\msk-r_ih\mse_i)]. \end{array}\]
For $a_{12}<0$ the form is defined by:
\begin{equation}\label{exp5.3}\begin{array}{ll} \displaystyle
 a_n^{(-)}(v,u) \ = \ \:\sum_{\mskd}& \displaystyle\left(
 \sum_{i = 1}^2\: a_{ii} \,v_{\mskd}(i,1)\,u_{\mskd}(i,1)
 +\:\sum_{i \ne j}\:a_{ij} v_{\mskd}(i,r_i)\,u_{\mskd}(j,r_j)\right .
 \\&+\left .  \displaystyle \sum_{i \ne j}\:a_{ij}\, \frac{r_i}{r_j}\:\big[
 v_{\mskd}(i,1) \,u_{\mskd}(i,1)  -v_{\mskd}(i,r_i)\,u_{\mskd}(i,r_i)\big]
 \right). \end{array}
\end{equation}
If the last summand of (\ref{exp5.3}) were omitted, then the resulting $A_n$ would
have neighborhoods ${\cal N}(\msx)$ containing seven grid-knots $\msx,\msx \pm
h\mse_i,\msx\pm h(r_1\mse_1-r_2\mse_2)$ as illustrated in Figure \ref{fig4.2}, as
well as the following four additional grid-knots, $\msx \pm hr_i\mse_i$. The second
line of (\ref{exp5.3}) causes cancellation of those matrix entries which would
correspond to four additional grid-knots.

For $a_{12}>0$ the form is obtained from the previous one by changing
$v_{\mskd}(i,\cdot), u_{\mskd}(i,\cdot), i=2$, into $\widehat{v}_{\mskd}(i,\cdot),
\widehat{u}_{\mskd}(i,\cdot), i=2$, respectively.

For the case of functions $\msx\mapsto a_{ij}(\msx)$ one naturally starts from the
just obtained expression for $a_n^{(\pm)}(v,u)$, changing the numbers $a_{ij}$ into
the numbers $a_{ij}(h\msk)$. In fact, the numbers $a_{ij}(h\msk+\mst)$, where $\mst$
are appropriately selected element of ${\bbR}^2$, are acceptable. The fastest
convergence of $h_n^2a_n(v,u)\to a(v,u)$ is a criterion which helps us to choose
the vectors $\mst$. It turns out that the best choice is $\mst(\msr)=(h/2)
(r_1\mse_1+r_2\mse_2)$ \cite{TS}, i.e. the values for which $h\msk+\mst
(\msr)$ are the mid-point of the rectangle $C_n(\msr,h\msk)$, with the lower
left vertex at $h\msk$ and the upper right vertex at $h\msk+r_1\mse_1+r_2\mse_2$.
Let us remark that the finite difference operators $\der_i(r_ih)u(h\msk)$ are
defined in terms of function values at the vertices of $C_n(\msr,h\msk)$. Thus, if
$a_{12}(\msx) \leq 0, \msx\in{\bbR}^2$, we get:
\[\begin{array}{ll}
 a_n(v,u) \ = \ \sum_{\msxd \in G_n}\:\Big(&\displaystyle \sum_{i = 1}^2\:
 a_{ii}(\msx+\mst({\bf 1}))\, v_{\mskd}(i,1)\,u_{\mskd}(i,1) \\  \displaystyle
 &+\:\sum_{i \ne j}\:a_{ij}(\msx+\mst(\msr))\, v_{\mskd}(i,r_i)\,u_{\mskd}(j,r_j)
 \\  &\displaystyle
 +\: \sum_{i \ne j}\:a_{ij}(\msx+\mst(\msr))\, \frac{r_i}{r_j}\: \big[
 v_{\mskd}(i,1) \,u_{\mskd}(i,1)  -v_{\mskd}(i,r_i)\,u_{\mskd}(i,r_i)\big]\Big),
 \end{array}\]
where $\msr =(r_1,r_2), {\bf 1}=(1,1)$.
By using the variational method, one obtains the entries $(A_n)_{\mskd \msld}, \msk,
\msl\in I_n$. At the present step of construction it is not necessary to write down
all the entries. In order to describe the influence of the parameters $r_1, r_2$
on the structure of entries we consider one group of entries:
\[ \big(A_n\big)_{\mskd \,\mskd+ \msed_1} \ = \  -\,\frac{1}{h^2}
 \:\Big[a_{11}(\msx+\mst({\bf 1}))\:-\: \frac{r_1}{r_2}\,
 \big|a_{12}(\msx+\mst(\msr))\big|\Big].\]
The matrices $A_n$ are of positive type iff the bracket has positive sign
for each $\msx\in G_n$. This is a condition on the functions $a_{ij}$, which
is implied by a particular choice of $r_1, r_2\in \bbN$.

For $a_{12}\geq0$ the form is obtained from the constructed one by changing
the following quantities. The finite differences $v_{\mskd}(i,\cdot),
u_{\mskd}(i,\cdot),i=2$, should be changed into $\widehat{v}_{\mskd}(i,\cdot),
\widehat{u}_{\mskd}(i,\cdot), i=2$, as in the case of constant coefficients, and
$\mst(\msm)$ should be changed into $\mss(\msm)=(h/2)(m_1\mse_1-m_2\mse_2)$.
With these changes, the obtained $A_n$ are symmetric matrices, possibly  of
positive type. We say ``possibly of positive type" since this depends on the
choice of $r_1, r_2$.

In the general case, the sign of $a_{12}$ is not constant on ${\bbR}^2$. Therefore,
we partition ${\bbR}^2$ into two subsets $\{\msx\in {\bbR}^2:a_{12}(\msx)
\leq 0 \}$ and $\{\msx\in {\bbR}^2:a_{12}(\msx)> 0 \}$. Then each of these
sets has to be partitioned further, where each of the partitioned classes is
characterized by a pair $r_1,r_2$, so that the resulting entries $(A_n)_{\mskd
\mskd\pm \msed_i}$ have negative values. The construction is carried out for a
class of functions $a_{ij}$ with moderate discontinuities.

\begin{assumption}\label{Ass4.1} Let there exist a finite index set $\EuScript{L}$,
a partition ${\bbR}^2 = \cup_{l\in\EuScript{L}} D_l$ and a diffusion tensor $a=
\{a_{ij}\}_{11}^{22}$
on ${\bbR}^2$ satisfying the strict ellipticity conditions (\ref{exp3.2}) and
the following additional discretization conditions:
\begin{description}\itemsep 0.cm
 \item{a)} There exist $\msq \in {\bbN}^2$, $n_0 \in {\bbN}$ and the corresponding
$h_0=2^{-n_0}$ so that each set $D_l$ is a connected union of cubes of the form
$\prod_{i=1}^d [x_i,x_i+q_ih_0)$. The matrix-valued function
$\msx \mapsto a(\msx)$ is continuous on ${\rm cl}(D_l)$ and either $a_{12}\geq 0$ or
$a_{12}\leq 0$ on $D_l$.
 \item{b)} To each $D_l$ there is associated a parameter $\msr(l) \in {\bbN}^2$,
such that the following inequality is valid:
\end{description}
\begin{equation}\label{exp4.8}\begin{array}{ll} \displaystyle
 \omega(a) \ = & \displaystyle\inf_n\: \min_{l\in {\cal L}}\, \min_{i=1,2}
 \, \inf \Big\{\, \inf_{\mszd \in S_n(\msrd(l),\msxd)\cap D_l}\,a_{ii}(\msz)\\ &
 \displaystyle -\: \frac{r_i(l)}{r_{m(i)}(l)}\, \sup_{\mszd \in S_n(\msrd(l),
 \msxd)\cap D_l} \,|a_{im(i)}(\msz)| \::\: \msx \in G_n \Big\} \ > 0.
 \end{array}
\end{equation}
where $m(i)=3-i$ and the rectangles $S_n(\msr(l), \msx)$ are defined by
(\ref{exp4.s}).
\end{assumption}

Condition b) is crucial in our construction of discretizations $A_n$ with the
structure of matrices of positive type.

The set $\EuScript{L}$ of Assumption \ref{Ass4.1} is partitioned into the subsets
$\EuScript{L}_{\mp}$, where $l \in \EuScript{L}_-$ means that $a_{12} \leq 0$
on $D_l$, and $l \in \EuScript{L}_+$ means that $a_{12} \geq 0$ on $D_l$.
It is convenient to use a representation
$a_n(v,u) = a_n^{(-)}(v,u) + a_n^{(+)}(v,u)$, where the form $a_n^{(-)}(v,u)$
contains the sums over the grid-knots in $D_l, l \in \EuScript{L}_-$ and
$a_n^{(+)}(v,u)$ contains the sums over the grid-knots in
$D_l, l\in\EuScript{L}_{+}$. Let us define
\begin{equation}\label{exp4.6}\begin{array}{l}
 \mst^{(\pm +)}(\msr) \ = \  \frac{h}{2}\,\big (\pm \,r_1(l)\mse_1 \,
 +\,r_2(l)\mse_2\big) \ \in S_n(\msr,{\bf 0}),\\
 \mst^{(\pm -)}(\msr) \ = \  \frac{h}{2}\,\big (\pm \,r_1(l)\mse_1 \,
 -\,r_2(l)\mse_2\big) \ \in S_n(\msr,{\bf 0}).
 \end{array}
\end{equation}
where $\msr=(r_1,r_2)$. The form $a_n^{(-)}(\cdot,\cdot)$ is defined by:
\begin{equation}\label{exp4.9}\begin{array}{ll}
 a_n^{(-)}(v,u)  = \ \sum_{l \in \EuScript{L}_-}\:&\displaystyle\Big\{
 \sum_{\msxd \in G_n(l)}\:\Big(\sum_{i = 1}^2\:
 a_{ii}(l,\msx+\mst^{(++)}({\bf 1}))\, \big(\der_i(h)v \big)(\msx) \big(
 \der_i(h)u \big)(\msx)  \\  &\displaystyle
 +\:\sum_{i,j=1,2,i \ne j}\:a_{ij}(l,\msx+\mst^{(++)}(\msr))\, \big(\der_i
 (r_i(l)h)v \big)(\msx)\big(\der_j(r_j(l)h)u \big)(\msx)  \\  &\displaystyle
 +\: \sum_{i,j=1,2,i \ne j}\:a_{ij}(l,\msx+\mst^{(++)}(\msr))\, \frac{r_i(l)}
 {r_j(l)}\: \Big[ \big(\der_i(h)v \big)(\msx) \big(\der_i(h)u \big)(\msx)\\
 &\:-\:\big(\der_i(r_i(l)h)v \big)(\msx) \big(\der_i(r_i(l)h)u
 \big)(\msx)\Big] \ \Big) \Big\},\end{array}
\end{equation}
where, as already noted, $\msr=(r_1,r_2)$, so that ${\bf 1}=(1,1)$.
The form $a_n^{(+)}(v,u)$ is obtained from $a_n^{(-)}(v,u)$ formally by replacing
$\EuScript{L}_-$ with $\EuScript{L}_+$, then $\der(mh)f$ with $\widehat{\der}(mh)f$
for each $m\in\{1,r_2\}$ and $f=v,u$, and $\mst^{(++)}$ with $\mst^{(+-)}$. Observe
that the forms $a_n^{(\mp)}(v,u)$ are the second degree polynomials in the
quantities $\der_i(qh), \widehat{\der}_i(qh)$ with $q\in\{1,r_1,r_2\}$.

We say that $\msx \in G_n$ is an {\em internal grid-knot} if
${\cal N}(\msx)\subset G_n\cap D_l$ for some $l$. All the other grid-knots are
called {\em boundary grid-knots}. For an internal grid-knot $\msx$ the acceptable
expressions for $(A_n)_{\mskd \msld}$ follow directly from the definition of
corresponding discrete forms ${\bf v}, {\bf u} \mapsto \lev {\bf v}|A_n{\bf u}\des$.
For a boundary grid-knot $\msx=h\msk$ the obtained expressions are complex, and
the calculated $(A_n)_{\mskd \msld}$ could break down the structure of matrices of
positive type. Therefore, one seeks simpler procedures for constructing
entries $(A_n)_{\mskd \msld}$ for boundary grid-knots $\msx=h\msk$. The results
of such a construction must be matrices $A_n$ which determine MJPs and the
convergence (\ref{exp1.5}) of MJPs determined by $A_n$ should also be ensured.

In order to write down the entries of $A_n$ corresponding to internal grid-knots,
we use (\ref{exp4.6}) and the following abbreviations:
\[\begin{array}{lll}
 \msw^{(\pm)}(l) &=& r_1(l)\mse_1 \,\pm\,r_2(l)\mse_2 \ \in I_n,\\
 a_{ij}^{(\alpha \beta)}(\msr) &=& a_{ij}( \msx + \mst^{(\alpha \beta)}(\msr)),
 \quad \alpha, \beta \in \{+,-\},\\
 \check{a}_{12}^{(-+)}(\msr) &=& a_{12}(\msx+\mst^{(++)}(\msr)-
 h\mse_1),\\
 \check{a}_{12}^{(+-)}(\msr) &=& a_{12}(\msx+\mst^{(++)}(\msr)-
 h\mse_2),\\
 \check{a}_{ii}^{(++)}(\msr) &=& a_{ii}^{(++)}(\msr),\quad
 \check{a}_{ii}^{(--)}(\msr) \ = \ a_{ii}^{(--)}(\msr).
\end{array}\]
When we apply variational method to the form $a_n^{(-)}$ defined by (\ref{exp4.9})
and the corresponding form $a_n^{(+)}$, we obtain entries of $A_n$. Thus, for an
internal grid-knot $h\msk$ the nontrivial off-diagonal entries of $A_n$ are:
\begin{equation}\label{exp4.23}\begin{array}{ll}\displaystyle
 \big(A_n\big)_{\mskd \mskd\pm \msed_1}& \displaystyle
 = \  -\,\frac{1}{h^2}
 \: \left\{\begin{array}{lll} \displaystyle
 a_{11}^{(\pm +)}({\bf 1})\:-\: \frac{r_1(l)}{r_2(l)}\,
 \big|\check{a}_{12}^{(\pm +)}(\msr)\big| &{\rm for}& a_{12} \leq 0, \
 {\rm on} \ D_l, \\
 \displaystyle
 a_{11}^{(\pm -)}({\bf 1})\:-\: \frac{r_1(l)}{r_2(l)}\,
 \big|\check{a}_{12}^{(\pm -)}(\msr)\big| &{\rm for}& a_{12} \geq 0, \
 {\rm on} \ D_l,  \end{array}
 \right . \\ \displaystyle
 \big(A_n\big)_{\mskd \mskd\pm \msed_2} & \displaystyle
 = \  -\,\frac{1}{h^2}
 \: \Big[ a_{22}^{(+\pm)}({\bf 1})\:-\: \frac{r_2(l)}{r_1(l)}\,
 \big|\check{a}_{12}^{(+\pm)}(\msr)\big|\Big].  \end{array}
\end{equation}
The entries corresponding to the grid-knots in the plane spanned by
$\mse_1, \mse_2$ have the structure:
\begin{equation}\label{exp4.24}\begin{array}{l}
 \begin{array}{l}\displaystyle
 \big(A_n\big)_{\mskd \mskd + \mswd^{(-)}(l)} \ = \ -\,\frac{1}{h^2r_1(l)r_2(l)}\:
 \big|a_{12}^{(+-)}(\msr)\big| \\ \displaystyle
 \big(A_n\big)_{\mskd \mskd - \mswd^{(-)}(l)} \ = \ -\,\frac{1}{h^2r_1(l)r_2(l)}\:
 \big|a_{12}^{(-+)}(\msr)\big|
 \end{array} \quad {\rm for}\quad a_{12} \leq 0 \ {\rm on} \ D_l, \\
 \begin{array}{l}\displaystyle
 \big(A_n\big)_{\mskd \mskd + \mswd^{(+)}(l)} \ = \ -\,\frac{1}{h^2r_1(l)r_2(l)}\:
 \big|a_{12}^{(++)}(\msr)\big| \\ \displaystyle
 \big(A_n\big)_{\mskd \mskd - \mswd^{(+)}(l)} \ = \ -\,\frac{1}{h^2r_1(l)r_2(l)}\:
 \big|a_{12}^{(--)}(\msr)\big| \\
 \end{array} \quad {\rm for} \quad a_{12} \geq 0 \ {\rm on} \ D_l, \end{array}
\end{equation}

If the quantities $\check{a}_{ij}^{(\alpha \beta)}(\msr), \alpha, \beta\in \{+,-\}$
in (\ref{exp4.23}) are replaced with $a_{ij}^{(\alpha \beta)}(\msr)$, the
convergence is still preserved. However, the quantities $a_{ij}^{(\alpha \beta)}
(\msr), \check{a}_{ij}^{(\alpha \beta)}(\msr)$ should
not be replaced with $a_{ij}(h\msk)$ since the resulting $(A_n)_{\mskd \msld}$
would be discretizations of $-\sum_{ij}a_{ij}\partial_i \partial_j$ and not
of $-\sum_{ij}\partial_i a_{ij}\partial_j$. This assertion can be easily proved
for the case of dimension $d=1$ and the diffusion tensor $a$, that
is twice continuously differentiable. We intend to compare the expressions
$A(x)u(x)=-(a(x)u(x)')'$ and $A^{cs}(x)u(x)=-a(x)u''(x)$ and their discretizations
on grids $G_n=\{hk:k\in \bbZ\}\subset \bbR$. The discretizations of $A^{cs}(x)u(x)$
are given by Lemma \ref{lem2.1}:
\[ \big(A^{cs}\big)_{kk\pm 1} \, =\, -h^{-2}a(hk), \qquad \big(A^{cs}\big)_{kk}
 \,=\, 2h^{-2}a(hk),\]
and consequently:
\[ \big((A^{cs})_n{\bf u}_n\big)_k \ = \ a(hk)\,\frac{2u(hk)-u(hk+h)-u(hk-h)}{h^2}.\]

For the discretizations of $A(x)u(x)$ we consider (\ref{exp4.23}). Let us represent
$a(x\pm h/2)$ by its Taylor polynomial of the second degree,
$a(x\pm h/2)=a(x)\pm(h/2)a'(x)+(h/2)^2a''(x)+r(\pm h,x)$, where the remainder
$r(\pm h,x)$ has the property $\lim_{h\to 0}h^{-2}r(\pm h,x)=0$. Therefore we get:
\[ \big(A\big)_{kk\pm 1} \, =\,-\,\Big(1+\frac{h^2}{4}\frac{a''(hk)}{a(hk)}
 \Big)\,\frac{1}{h^{2}}\,a(hk)\, \mp \frac{1}{2h}a'(hk)+o(\pm h,x), \]
where $\lim_{h\to 0}o(\pm h,x)=0$. Hence:
\[\begin{array}{lll} \big(A_n{\bf u}_n\big)_k &=& \displaystyle
 \Big(1+\frac{h^2}{4}\frac{a''(hk)}{a(hk)}
 \Big)\,a(hk)\, \frac{2u(hk)-u(hk+h)-u(hk-h)}{h^2}\\
 &+& \displaystyle a'(hk)\frac{u(hk+h)- u(hk-h)}{2h}\,+\,\tilde{o}(h,x).
 \end{array}\]
From the obtained expressions we have:
\[ \big(A_n{\bf u}_n\big)_k\,-\,
\Big(1+\frac{h^2}{4}\frac{a''(hk)}{a(hk)}\Big)\big((A^{cs})_n{\bf u}_n\big)_k
 \ = \  \big(B_n{\bf u}_n\big)_k\,+\,\tilde{o} (h,x),\]
where $\big(B_n{\bf u}_n\big)_k$ are the discretizations of $-a'(x)u'(x)$. In other
words, we have obtained discretizations of the expression $A(x)=A^{cs}(x)-
a'(x)(d/dx)$ as we should have.

Now we can describe the construction which gives a satisfactory result for both
types of grid-knots, internal and boundary ones. For each $\msx=h\msk \in
{\rm cls}(D_l)$ the entries $(A_n)_{\mskd \msld}$ are constructed by the rules
(\ref{exp4.23}), (\ref{exp4.24}). If ${\cal N}(\msx)\subset {\rm cls}(D_l)$ there
is nothing more to adjust. If ${\cal N}(\msx)\cap{\rm cls}(D_l) \notsubseteq
{\cal N}(\msx)$ then the quantities $a_{12}^{(\alpha \beta)}(\msr),
\check{a}_{12}^{(\alpha \beta)}(\msr)$, where $\alpha, \beta \in \{+,-\}$,
should be replaced by zeros in all the expressions (\ref{exp4.23}),
(\ref{exp4.24}). Let us point out that otherwise the entries (\ref{exp4.23})
for the case of ${\cal N}(\msx)\cap{\rm cls}(D_l) \notsubseteq {\cal N}(\msx)$
could have a wrong sign, i.e. it could happen that
$(A_n)_{\mskd \msld}> 0$ for some $\msk \ne \msl$. In order to
avoid such undesired features, one therefore omits the terms of entries in
the expressions (\ref{exp4.23}), (\ref{exp4.24}), which are proportional
to $|a_{12}|$. This adjustment procedure is equivalent to the assumption that the
function $a_{12}$ is zero in a neighborhood of the boundary
$\Gamma=\cup_l \partial D_l$.

This determines the rules of construction of discretizations of a
generalized diffusion for which the diffusion tensor satisfies Assumption
\ref{Ass4.1}. The above described construction of entries $(A_n)_{\mskd \msld}$
at boundary grid-knots can be justified as follows.

A numerical neighborhood ${\cal N}(\msx)$, where $\msx \in D_l \cap G_n, l\in
{\cal L}$ is contained in the closed rectangle
\begin{equation}\label{exp3.1a}
 S_n(\msr(l),\msx) \ = \ \prod_{i=1}^d \: [x_i - h_nr_i(l)\,,\,x_i+h_nr_i(l)].
\end{equation}
The union of all such rectangles centered at boundary grid-knots $\msx \in \Gamma
\cap G_n$ is denoted by $S_n(\Gamma)$. This is a closed set covering $\Gamma$.
Because of $G_n\subset G_{n+1}$ we have $S_{n+1}(\Gamma)\subset S_n(\Gamma)$
and the identity $\Gamma =\cap_nS_n(\Gamma)$. For each $n \in \bbN$, we approximate
the original diffusion tensor $a$ by the tensor $a^{(n)}$ defined as follows:
\begin{equation}\label{exp4.10}\begin{array}{lllll}
 a_{ij}^{(n)}(\msx) &=& a_{ij}(\msx) &{\rm for} &\msx \in {\bbR}^2\setminus
 S_n(\Gamma),\quad i,j \in \{1,2\},\\
 a_{ii}^{(n)}(\msx) &=& a_{ii}(\msx) &{\rm for} &\msx \in S_n(\Gamma),
 \quad i \in \{1,2\},\\
 a_{12}^{(n)}(\msx) &=& 0            &{\rm for} &\msx \in S_n(\Gamma).\end{array}
\end{equation}
The defined diffusion tensors $a^{(n)}$ determine differential operators
$A_0^{(n)}(\msx)$ which approximate the differential operator $A_0(\msx)$. Each
$A_0^{(n)}(\msx)$ determines a generalized diffusion $X^{(n)}(\cdot)$.
The generalized diffusions $X^{(n)}(\cdot)$ converge in distribution to the original
diffusion $X(\cdot)$. Therefore it suffices to consider the discretizations
$A_n^{(n)}$ which are defined in terms od $a_{ij}^{(n)}$ as described above.
The convergence in distribution of diffusions $X^{(n)}(\cdot)$ to $X(\cdot)$
follows from the convergence of the corresponding semigroups $U^{(n)}(\cdot)$
to the semigroup $U(\cdot)$ in $\dot{C}({\bbR}^d)$. To see this,
we first apply Corollary \ref{Cor3.1} in order to prove the uniform
boundedness $\sup_{t\leq 1}\Vert U^{(n)}(t)\Vert^{(\alpha)}\leq
\beta$, where $\alpha, \beta$ do not depend on $n$. Then we use the
standard methods of Sobolev spaces in order to prove the convergence of the
semigroups in $L_2({\bbR}^d)$. The two above properties are combined in order
to prove the convergence of semigroups in $\dot{C}({\bbR}^d)$ as in the last
step of proof of Theorem \ref{Th5.2}.

\subsection{Construction for $d \geq 3$}
The goal of the overall analysis is to find those discretizations $A_n$ of the
differential operator $A(\msx)$ which have the structure of matrices of positive
type. Here we describe a general approach, which is based on reduction to a
finite number of two-dimensional problems.

The index set of pairs $I(d) = \{ \{ij\} : i < j,\:
i,j = 1,2,\ldots,d, i \ne j\}$ has the cardinal number $m(d)=d(d-1)/2$.
To each index $\{kl\} \in I(d)$ we associate three coefficients,
\begin{equation}\label{exp4.12}
 a_{kk}^{\{kl\}} = \frac{1}{d-1}\, a_{kk}, \quad a_{ll}^{\{kl\}} =
 \frac{1}{d-1}\, a_{ll}, \quad a_{kl}^{\{kl\}} = a_{kl},
\end{equation}
and a two-dimensional bilinear form $a^{\{kl\}}(\cdot, \cdot)$,
\[ a^{\{kl\}}(v,u) \ = \sum_{i,j \in \{k,l\}}\: \int_{{\bbR}^d} \: a_{ij}^{\{kl\}}
 (\msx) \: \partial_i v(\msx)\, \partial_j u(\msx) \,d\msx.\]
Clearly, for each pair $v, u \in C_0^{(1)}({\bbR}^d)$
the following identity is valid:
\begin{equation}\label{exp4.13}
 a(v,u) \ = \ \sum_{\{kl\} \in I(d)} \, a^{\{kl\}}(v,u).
\end{equation}
To each of the forms $a^{\{kl\}}(\cdot, \cdot)$ we associate a sequence of
forms $a_n^{\{kl\}}(\cdot,\cdot)$ and matrices $A_n^{\{kl\}}$ constructed by using
schemes in two dimensions from the previous subsection. Then the matrix
\begin{equation}\label{exp4.14}
 A_n \ = \ \sum_{\{kl\} \in I} \:A_n^{\{kl\}} ,
\end{equation}
is a discretization of $A_0(\msx)$. If each $A_n^{\{kl\}}$ has the structure of a
matrix of positive type, then $A_n$ is also a matrix of positive type.
However, $A_n$ can have the structure of a matrix of positive type, even though
no $A_n^{\{kl\}}$ is a matrix of positive type. This important property, which
enables us to construct matrices $A_n$ of positive type for the case $d \geq 3$,
can be proved from the structure of entries $(A_n)_{\mskd \mskd\pm \msed_i}$.

First we consider the entry $(A_n)_{\mskd \mskd+\msed_1}$ defined by
(\ref{exp4.14}) in the case that $a_{ij}<0$ for all $i\ne j, \ i,j=1,2,\ldots,d$.
In addition, in order
to write down expressions as simply as possible, the index $l\in {\cal L}$
is omitted from the notations. The contribution from the sum of entries
$\big(A_n^{\{kl\}}\big){\mskd \mskd+\msed_1}$ to the entry $(A_n)_{\mskd
\mskd+\msed_1}$ has the following form:
\[ \frac{1}{d-1}\sum_{s\geq 2}a_{11}\Big(\msx+\frac{h}{2}(\mse_1+\mse_s)
 \Big) \ - \ {\rm terms~containing~~} a_{12}, a_{13},\ldots, a_{1d}.\]
Similarly we can describe the terms containing $a_{ii}$ for any $i=1,2,\ldots,d$, and any
$l\in {\cal L}$. The corresponding sum contributing to $(A_n)_{\mskd \mskd+\msed_i}$
has the following general form:
\begin{equation}\label{exp4.7}
 \omega_n(a_{ii},\msx) \ = \ \frac{1}{d-1}\:\sum_{s=1,s \ne i}^d \:
 a_{ii}(h\msk+h\msm_{ii}(l,s)),
\end{equation}
where $\msm_{ii}(l,s)$ are defined by the rules of construction of (\ref{exp4.14}),
(\ref{exp4.23}) and (\ref{exp4.24}). The terms proportional to $a_{is}, s\ne i$,
are summed with the just defined $\omega_n(a_{ii},\msx)$ as shown in the next
description of the obtained results.

\begin{procedure}\label{dsp4.1} Let there be given a partition ${\bbR}^d=\cup_l D_l$
into a finite number of subsets $D_l$ such that all the functions $a_{ij}$ are
uniformly continuous on each $D_l$, and the functions $a_{ij}, i\ne j, i,j=
1,2,\ldots,d$, do not change sign on $D_l$. Let a parameter $\msr(l)=(r_1(l), r_2(l),
\ldots, r_d(l)) \in {\bbN}^d$ be assigned to each $D_l$ and the matrices
$A_n$ on $G_n$ be constructed by the rules (\ref{exp4.23}), (\ref{exp4.24})
and (\ref{exp4.14}). Then their entries have the following properties:
\begin{enumerate}
 \item Entries of $(A_n)_{\mskd \msld}, \msx = h \msk$, $\msk \in {\bbZ}^d$,
$h\msl \in {\cal N}(\msx)$, are linear combinations of
$a_{ij}(\msx_{ij}(n,\msx,l))$ where $\msx_{ij}(n,\msx,l) = h\msk + h\msm_{ij}(l,s)$,
and where the elements $\msm_{ij}(l,s)\in {\bbR}^d$ for $i,j,s =1,2,\ldots,d,
l \in \EuScript{L}$, do not depend on $n$.
 \item For each grid-knot $\msx = h \msk$: $(A_n)_{\mskd \mskd} = -
\sum_{\msld}  (A_n)_{\mskd \msld}$.
 \item For each $\msx = h \msk \in {\rm cls}(D_l)$ entries in the coordinate
directions $\msx \pm h\mse_i$ are defined by:
\[ \big(A_n \big)_{\mskd \mskd\pm \msed_i} \:=\:
 -\,\frac{1}{h^2}  \,\Big [ \omega_n(a_{ii},\msx)
 -\,\sum_{m=1,m \ne i}^d\:\frac{r_i(l)}{r_m(l)}\,
 |a_{im}(\msx_{im}(n,\msx,l))|\Big]. \]
 \item For each $l \in \EuScript{L}$ the entries of $A_n$ in the plane spanned
by $\mse_i,\mse_j$ are defined by using elements
$\msz_{ij}(l) = r_i(l)\mse_i - r_j(l)\mse_j \in {\bbZ}^d$ (if $a_{ij}\leq 0$ on
$D_l$) or elements $\msz_{ij}(l) = r_i(l)\mse_i + r_j(l)\mse_j \in {\bbZ}^d$, (if
$a_{ij}\geq 0$ on $D_l$), as follows:
\[ \big(A_n\big)_{\mskd \mskd \pm \mszd_{ij}(l)} \ = \ -\:\frac{1}{h^2
 r_i(l)  r_j(l)} \: |a_{ij}(\msx_{ij}(n,\msx,l))|.\]
\end{enumerate}
\end{procedure}

An appropriate choice of the parameters $\msr(l)$ follows from Theorem \ref{Th4.1}.

Some special features regarding the structure of the sets ${\cal N}(\msx), \msx =
h\msk \in G_n$, should be pointed out. If $a_{ij}\ne 0, i \ne j$, then the maximal
number of elements in ${\cal N}(\msx)$ is $1+d+d^2$. In this case the set
${\cal N}(\msx)$ consists of its center, $2d$-grid-knots in the coordinate
directions, and 2 grid-knots in each of $d(d-1)$ two-dimensional
planes. Since there can be at most two grid-knots in a two-dimensional plane
(Property 4. of Discretization procedure \ref{dsp4.1}),
the entries of $A_n^{(rs)}$ have the following property. Let the pairs $\mse_r,
\mse_s$ and $\mse_s, \mse_t$ span
two-dimensional planes and let $A_n^{(rs)}, A_n^{(st)}$ be the corresponding
discretizations which are constructed by using parameters $\msr^{(rs)},
\msr^{(st)}$. Then for the construction defined by Discretization procedure
\ref{dsp4.1} the following identity $(\msr^{(rs)})_s =(\msr^{(st)})_s$ must be
valid.

\subsection{Lower order differential operators}
The discretizations of differential operator $B(\msx)=\msb(\msx)\nabla$
are denoted by $B_n$. The following general rule should be obeyed.
A positive diagonal entry and a non-positive off diagonal entry is
associated to each $\msx \in G_n$ for which $\msb(\msx) \ne {\bf 0}$.
Let us define the sets $\EuScript{K}(i,-)=\{\msx \in G_n: b_i(\msx) < 0\}$
and analogously $\EuScript{K}(i,+)=\{\msx \in G_n: b_i(\msx) > 0\}$.
Then the discretizations of $(v|Bu)$ are defined by
\[ b_n(v,u) \ = \ \sum_i\:\Big[\sum_{\msxd\in \EuScript{K}(i,-)}\:
 b_i(\msx)\,v(\msx)\,\big(\der_i(h)u\big)(\msx)\,+\,
 \sum_{\msxd\in \EuScript{K}(i,+)}\:
 b_i(\msx)\,v(\msx)\,\big(\der_i(-h)u\big)(\msx)\Big].\]
These forms have to be summed with the forms (\ref{exp4.13}) in order to get
discretizations of the original form (\ref{exp3.4}). If discretizations
$(A_0)_n$ of $A_0(\msx)$ have the structure of matrices of positive type, then
obviously $(A_0)_n+B_n$ maintain this structure. The so defined discretizations
of $B$ are usually called {\em upwind schemes}.

The constructed forms $a_n$ of this section are discretizations of the form
(\ref{exp3.4}). At the present level of analysis the constructed discretizations
can be justified by the limit $a(v,u) = \lim_n h^d a_n(v,u)$, being valid for
any pair $v,u \in C_0^{(1)}({\bbR}^d)$.

\subsection{Summarized results of the construction}
\begin{theorem}\label{Th4.1} Let there be given a partition ${\bbR}^d=\cup_l D_l$
into a finite number of connected sets $D_l$, each being the union of cubes
$C_m(\msp,\msx)$ with some fixed $m$, so that the functions $a_{ij}$ fulfil
(\ref{exp3.2}) and the following additional conditions:
\begin{description}\itemsep -0.1cm
 \item {a)} The functions $a_{ij}$ are uniformly continuous on $D_l$
and $a_{ij}, i\ne j$, do not change sign on $D_l$.
 \item {b)} For each pair $i,j$ the limit $\lim_{|\msxd|\to \infty}a_{ij}
(\msx)$ has a constant value.
 \item {c)} The matrix-valued function $\msx\mapsto \hat{a}(\msx)$, defined by
(\ref{exp1.2}), is strictly positive definite on ${\bbR}^d$, i.e. $(\msz|\hat{a}
(\msx)\msz)\geq \beta |\msz|^2$ for some $\beta >0$ and all $\msz,\msx\in {\bbR}^d$.
\end{description}
Then there exist discretizations $A_n$ which are constructed by the rules of
Discretization procedure \ref{dsp4.1}, such that $A_n$ are matrices of positive type.
\end{theorem}

{\Proof} If for each $D_l$ we choose the parameters $r_i(l)$ so that the
entries $\big(A_n \big)_{\mskd \mskd\pm \msed_i}$ of item 3. of Discretization
procedure \ref{dsp4.1} have all negative values, then the rules of construction
(\ref{exp4.24}) ensure the existence of $A_n$ with the structure of matrices
of positive type. It remains to justify the existence of such
parameters $\msr(l)\in {\bbN}^d$ for each $l\in {\cal L}$.
Let us consider a set ${\rm cls}(D_l)$ and the quantity:
\begin{equation}\label{exp4.15}\begin{array}{ll} \displaystyle
 \omega(a) \ = & \displaystyle\inf_n\: \min_{l\in {\cal L}}\, \min_{i=1,2,\ldots,d}
 \, \inf \Big\{\, \inf_{\mszd \in S_n(\msrd(l),\msxd)\cap D_l}\,a_{ii}(\msz)\\ &
 \displaystyle -\: \sum_{m=1, m\ne i}^d\,
 \frac{r_i(l)}{r_{m}(l)}\, \sup_{\mszd \in S_n(\msrd(l),
 \msxd)\cap D_l} \,|a_{im}(\msz)| \::\: \msx \in G_n \Big\}.\end{array}
\end{equation}
If $\omega(a)>0$, the chosen parameters $\msr(l)$ ensure the positive
value of the brackets in item 3. of Discretization procedure \ref{dsp4.1}.
If $\omega(a)\leq0$, the partition should be refined until the condition
$\omega(a)>0$ is achieved. In accordance with Lemma \ref{lem2.1}, for each $\msx
\in G_n$ there exist $(r_1(\msx), r_2(\msx),\ldots,r_d(\msx))$ such that
$a_{ii}(\msx)-\sum_{m\ne i}(r_i(\msx)/r_m(\msx))|a_{im}(\msx)|>0$.
Due to the uniform continuity of the functions $a_{ij}$
on the sets ${\rm cls}(D_l)$, and b), the described procedure results with a
desirable result after a finite number of steps. {\QED}

The uniform continuity of functions $a_{ij}$ on $D_l$, and the inequality
$\omega(a)>0$, where $\omega(a)$ is defined by (\ref{exp4.15}), imply
another important property of matrices $A_n$ that are
constructed by the above described procedure. There exists $\sigma_0>0$,
independent of $n$, such that
\begin{equation}\label{exp72.10}
 \Big |\big(A_n\big)_{\mskd \mskd\pm \msed_i}\Big| \ \geq \ \frac{\sigma_0}{h^2}
 , \quad \msx =h\msk \in G_n.
\end{equation}
The described construction of matrices $A_n$ for which (\ref{exp72.10}) holds is
called {\em admissible method}\index{admissible method}, anticipating that the
obtained $A_n$ have all the necessary properties for the convergence
of corresponding MJPs to generalized diffusion. Let us recall $U_i=U_i(1)$.

\begin{lemma}\label{lem72.1} Let the matrices $A_n$ on $G_n$ be discretizations
of $A_0 = -\sum \partial_i a_{ij}\partial_j$ by an admissible method. Then the
matrices $A_n$ are irreducible. If in addition, the matrices $A_n$ are symmetric,
then there exist positive numbers $\Md,\Mg$ such that the following inequalities
\[ \Md\,\sum_i \lnorm U_i {\bf u}\lnorm_2^2 \ \leq \
 \lev {\bf u}\,|\,A_n{\bf u}\des
 \ \leq \ \Mg\,\sum_i \lnorm U_i {\bf u}\lnorm_2^2,\]
are valid uniformly with respect to $n\in \bbN$.
\end{lemma}

{\Proof} It is easy to check that a symmetric matrix $A_n$ of positive type is
positive semidefinite, i.e. $\lev {\bf u} |A_n {\bf u}\des\geq 0$.
Let us consider the tensor valued functions $a, b$ where $a(\msx) = \{a_{ij}(\msx)
\}_{11}^{dd}$ and $b(\msx)$ is defined by $b_{ij}(\msx)=a_{ij}(\msx)-\kappa
\delta_{ij}$. The corresponding auxiliary tensors of (\ref{exp1.2}) are denoted by
$\hat{a}, \hat{b}$. Due to the strict positive definitness of the matrices $a,
\hat{a}$ one can choose $\kappa > 0$
sufficiently small so that $b, \hat{b}$ are also positive definite on ${\bbR}^d$.
Let us define matrices $H_n$ on $G_n$ by the following non-trivial entries:
\[ \big(H_n\big)_{\mskd \mskd} \ = \ \frac{2d }{h^2}, \quad
 \big( H_n \big)_{\mskd \mskd \pm \msed_i} \ = \ -\frac{1}{h^2}, \
 i = 1,2,\ldots,d.\]
We have $\lev {\bf u}|H_n{\bf u}\des = \sum_{i=1}^d \lnorm U_i{\bf u}\lnorm_{2}^2$.
In accordance with the construction of tensors $b, \hat{b}$ and the inequality
(\ref{exp72.10}), the matrices $B_n = A_n - \kappa H_n$ are also of positive type.
Since the symmetric matrix $B_n$ of positive type is necessarily positive
semidefinite, i.e. $\lev {\bf u} |B_n {\bf u}\des \geq 0$ for any ${\bf u}
\in l_0(G_n)$, we have:
\[\lev {\bf u} \,|\,A_n \,{\bf u}\des \ = \ \kappa\lev {\bf u} \,|\,H_n \,{\bf u}
 \des \:+\: \lev {\bf u} \,|\,B_n \,{\bf u}\des \
 \geq \ \kappa \sum_{i=1}^d \lnorm U_i{\bf u}\lnorm_{2}^2, \]
proving the lower bound of the assertion. The upper bound follows from
(\ref{exp4.9}) and (\ref{exp4.13}). The irreducibility follows from the graph
theory, since any two grid-knots $\msx_0,\msy\in G_n$ can be connected by a path
of the form $\{\msx_0, \msx_1,\ldots,\msx_m,\msy\}\subset G_n$ such that
${\cal N}(\msx_{k-1})\cap{\cal N}(\msx_k) \ne \emptyset$ for $k=1,2,\ldots,m$. {\QED}

For a differential operator $A(\msx)=A_0(\msx)+\msb(\msx)\nabla$ with non-constant
coefficients, the functions $G_n\owns \msx\mapsto \delta(A,n,\msx)p(\msx)$ of
(\ref{exp4.0}) are not identically zero on $G_n$. In the proof of convergence
in Section \ref{sec5} the following weaker result is therefore used:

\begin{lemma}\label{lem5.1}
Let the differential operator $A(\msx)$ satisfy the conditions of Theorem \ref{Th4.1}
and let the following additional conditions be valid for some $\alpha \in (0,1)$:
\begin{description}\itemsep -0.15cm
 \item{a)} Functions $a_{ij}$ belong to the class $C^{(1+\alpha)}$ on ${\bbR}^d$.
 \item{b)} Functions $b_i$ belong to the class $C^{(\alpha)}$ on ${\bbR}^d$.
\end{description}
Then there exists a positive constant $\kappa$, depending
on $\Mg, \Vert \msb\Vert_\infty$, $\Vert a_{ij}\Vert^{(1+\alpha)}$ and
$\Vert b_{i}\Vert^{(\alpha)}$, but not on $n$, such that
\[ \sup_{\msxd\in G_n} \,|\delta(A,n,\msx)f(\msx)| \ \leq \ \kappa \, h^\alpha\,
 \Vert f\Vert^{(2+\alpha)}\]
for any $f\in C^{(2+\alpha)}({\bbR}^d)$.
\end{lemma}

This result can be easily checked by calculating $A_n{\bf f}_n$ directly.

Let us mention the following. If the functions $a_{ij}, b_i$ are uniformly of
the class $C^{(1+\alpha)}$ and $C^{(\alpha)}$ on $D_l$, respectively, and
$f\in C^{(2+\alpha)}$ then
$\delta(A,n,\msx)f(\msx)$ converges to zero on grid-knots $\msx\in {\rm int}(D_l)
\cap\big(\cup_n G_n\big)$. Otherwise, $n\mapsto \delta(A,n,\msx)f(\msx)$ is not
bounded as $n$ increases. Nevertheless, the convergence in $W_2^1$-spaces \cite{LR2}
is ensured as usually.

\section{CONVERGENCE OF MJPs}\label{sec5}
The convergence of MJPs to generalized diffusion is analyzed here in terms of
the criterion (\ref{exp1.5}). Therefore, we first need to define explicitly the
mappings $\Phi_n: \dot{l}_\infty(G_n) \mapsto \dot{C}({\bbR}^d)$.

An element (column) ${\bf u}_n \in l(G_n)$ can be associated to a continuous
function on ${\bbR}^d$ in various ways. In the current setting we define a mapping
$l(G_n) \mapsto C({\bbR}^d)$ in terms of {\em hat functions}. Let $\chi$ be the
canonical hat function on ${\bbR}$, centered at the origin and having the support
$[-1,1]$:
\[ \chi(x) \ = \ \left\{\begin{array}{lll} 1+x &{\rm for}& x\in [-1,0],\\
 1-x &{\rm for}& x\in [0,1],\\ 0 &{\rm for}& |x|>1.\end{array}\right .\]
Then $z \mapsto \phi(h,x,z) = \chi(h^{-1}(z-x))$ is the hat function on ${\bbR}$,
centered at $x \in {\bbR}$ with support $[x-h,x+h]$. The functions $\msz  \: \mapsto
\: \phi_{\mskd}(\msz) \:=\:\prod_{i=1}^d \:\phi(h,x_i,z_i),x_i = hk_i, i=1,2,
\ldots,d$, are the corresponding $d$-dimensional hat functions with support
$S_n({\bf 1},\msx) = \prod_i [x_i-h, x_i+h]$. The functions $\phi_{\mskd}(\cdot)
\in l(G_n)$, span a linear space, denoted by $E(n,{\bbR}^d)$.
Let ${\bf u}_n \in l(G_n)$ have the entries $u_{n \mskd} = ({\bf u}_n)_{\mskd}$.
Then the function $u(n) =\sum_{\mskd \in I_n} u_{n  \mskd} \phi_{\mskd}$
is an element of $E(n,{\bbR}^d)$ and defines an embedding of grid-functions into the
space of continuous functions. We denote the corresponding mapping by
$\Phi_n : l(G_n) \mapsto E(n,{\bbR}^d)$ and write
\begin{equation}\label{exp2.5}
 u(n) \ = \ \Phi_n\,{\bf u}_n \ = \ \sum_{\mskd}\:
 \big({\bf u}_n\big)_{\mskd}\:\phi_{\mskd}.
\end{equation}
The inverse mapping $\Phi_n^{-1} : E(n,{\bbR}^d)  \mapsto l(G_n)$ is defined by
$\Phi_n^{-1}\big(\sum u_{\mskd}\phi_{\mskd}\big)={\bf u}$, where the column
${\bf u}$ has the entries $u_{\mskd}$. It is obvious that the spaces $l(G_n)$ and
$E(n,{\bbR}^d)$ are isomorphic with respect to the pair of mappings
$\Phi_n, \Phi_n^{-1}$. Since $h_n=2^{-n}$, it is clear that $E(n,{\bbR}^d)
\subset E(n+1,{\bbR}^d)$ and that the space of functions $\cup_n E(n,{\bbR}^d)$
is dense in $L_p({\bbR}^d), p \in [1,\infty)$, as well as in $\dot{C}({\bbR}^d)$.
Let us mention that $\sum_{\mskd}\phi_{\mskd} = 1$ on ${\bbR}^d$. For two
functions $v(n), u(n)\in E(n,{\bbR}^d)$ we have $(v(n)|u(n)) = h^d \sum_{\mskd
\msld} s_{\mskd \msld} v_{\mskd} u_{\mskd}$ where $s_{\mskd \msld} = \Vert
\phi_{\mskd}\Vert_1^{-1}(\phi_{\mskd}|\phi_{\msld})$. The numbers $s_{\mskd \msld}$
do not depend on $n$ and the following identity is valid: $\sum_{\msld} s_{\mskd
\msld}= 1$. Thus we have $\Phi_n^{-1} \Phi_n=I$ in $l(G_n)$ and $\Phi_n\Phi_n^{-1}=I$
in $E(n,{\bbR}^d)$. Let $P(n)$ be the projector onto $E(n,{\bbR}^d)$ defined by
$f\mapsto P(n)f = \Phi_n {\bf f}_n$. The mapping $\Phi_n^{-1}$ can be extended
from $E(n,{\bbR}^d)$ to $\dot{C}({\bbR}^d)$ by defining $\Phi_n^{-1}f=\Phi_n^{-1}
P(n)f$. Thus we have $\Phi_n\Phi_n^{-1}=P(n)$ in $\dot{C}({\bbR}^d)$.

If $F_n$ is a matrix on $G_n$, then $F(n)=\Phi_n F_n\Phi_n^{-1}$ is a linear
operator in the linear space $E(n,{\bbR}^d)$. It is easy to verify that
$\Vert F(n)\Vert_p=\lnorm F_n\lnorm_p$ for $p = 1,\infty$, where the norm
$\Vert \cdot\Vert_p$ is induced by the restriction of $L_p({\bbR}^d)$
to $E(n,{\bbR}^d)$. By applying the interpolation Rietz-Thorin theorem \cite{BL}
we get $\Vert F(n)\Vert_p\leq\lnorm F_n\lnorm_p$ for $p \in [1,\infty]$.

The objective of the analysis in this section is the comparison of the Feller
semigroup $U(\cdot)$ in $\dot{C}({\bbR}^d)$ and the matrix semigroups $U(n,t)=\Phi_n
\exp(-A_n t)\Phi_n^{-1}$ in $\dot{E}(n,{\bbR}^d) =E(n,{\bbR}^d)\cap
\dot{C}({\bbR}^d)$, leading to a proof of (\ref{exp1.5}). For this purpose we
consider the following initial value problems (IVP):
\begin{equation}\label{exp5.2}\begin{array}{l}
 \big(\partial_t+A(\msx)\big) u(t,\msx) \ = \ 0, \quad u(0,\msx)=u_0(\msx),\\
 \dot{{\bf u}}_n(t)+ A_n{\bf u}_n(t) = {\bf 0}, \quad {\bf u}_n(0)=
 {\bf u}_{0n}, \quad n\in \bbN.\end{array}
\end{equation}
where ${\bf u}_{0n}$ are the discretizations of $u_0$.
By using the standard methods in the Sobolev space $W_2^1({\bbR}^d)$, we can first
prove a result which is weaker than (\ref{exp1.5}). Let $A(\msx)$ be defined by
(\ref{exp3.1}) and Assumption \ref{Ass3.1}, and let $A_n$ be its discretizations
on $G_n$, constructed by the rules of Section \ref{sec4}. We consider $u(t)=U(t)f$
for $f\in \dot{C}({\bbR}^d)\cap W_2^1({\bbR}^d)$, the function $f(n)$ defined by
(\ref{exp2.5}), and $u(n,t)=U(n,t)f(n)$. As proved in \cite{LR2}, for each
specified $f$ we have  $\lim_{n\to\infty}\:\Vert u(t)-u(n,t)\Vert_{2,1}=0$,
uniformly on segments of $[0,\infty)$.

The function $u(t)=U(t)f$ for $f\in \dot{C}({\bbR}^d)\cap W_2^1({\bbR}^d)$ is
continuous on ${\bbR}^d$, uniformly on segments of $[0,\infty)$, as shown in
\cite{LSU}. The corresponding functions $u(n,t)$ are also continuous on ${\bbR}^d$,
uniformly on segments of $[0,\infty)$, as follows from their structure, $u(n,t)=
\Phi_n^{-1}\exp(-A_nt){\bf f}_n$.
Due to the just described convergence in $W_2^1({\bbR}^d)$,
the sequence $\mathfrak{U}=\{u(n,t):n\in \bbN\}$ is bounded in
$W_2^1({\bbR}^d)$, uniformly on segments of $[0,\infty)$. Hence, this sequence
has a subsequence converging a.e. to $u$ on ${\bbR}^d$. However, we need to
show the uniform convergence. As is usual in such problems, the uniform
convergence can be proved eventually for an appropriately selected subsequence of
$\mathfrak{U}$. In our
approach, the outline of proof of (\ref{exp1.5}) is as follows. The original
differential operator $A(\msx)$ is approximated by differential operators
$A^{(m)}(\msx)$ with smooth coefficients. The corresponding semigroups are denoted
by $U^{(m)}(\cdot)$ in $\dot{C}({\bbR}^d)$ and $U^{(m)}(n,\cdot)$ in
$\dot{l}_\infty(G_n)$. The limit (\ref{exp1.5}) is then proved for each $m\in\bbN$.
Finally, by applying the diagonalization argument to the sequence $\{u^{(m)}(n,t):
m,n\in\bbN\}$ we get the desired result. The so outlined steps of the proof are
performed in the next two subsections.

\subsection*{Convergence for smooth coefficients}
Let the differential operator (\ref{exp3.1}) have the coefficients $b_i$ that belong
to the class $C^{(\alpha)}$ and the coefficients $a_{ij}$ that belong to the class
$C^{(1+\alpha)}$ so that it can be represented as,
\[ A(\msx) \ = \ -\,\sum_{ij=1}^d\,a_{ij}(\msx)\partial_i\partial_j \,+\,
 \sum_{i=1}^d\, b_i'(\msx)\partial_i,\]
where $b_i'(\msx)=b_i(\msx)- \sum_j \partial_j a_{ij}(\msx)$. Hence, $A(\msx)$ can
be represented as an elliptic operator in non-divergence form with coefficients
belonging to the class $C^{(\alpha)}$ on ${\bbR}^d$. This form of $A(\msx)$
makes it possible to use results on the existence of a strongly continuous
semigroup in $\dot{C}^{(2+\alpha)}({\bbR}^{d})=C^{(2+\alpha)}({\bbR}^{d})\cap
\dot{C}({\bbR}^{d})$ as developed by Solonnikov
\cite{So} (a detailed exposition of results can be found in \cite{LSU}, Sections 13
and 14 of Chapter 4). Thus we have $\Vert U(t)\Vert^{(2+\alpha)}
\leq \exp(\sigma t)$ for some $\sigma\geq 0$.

The matrices $A_n$ of Section \ref{sec4} approximate $A(\msx)$ as described
in Lemma \ref{lem5.1}. We have $\big(\Phi_n^{-1} A - A_n \Phi_n^{-1}\big)f(h\msk)=
\delta(A,n,h\msk)f(h\msk)$ so that:
\begin{equation}\label{exp5.5}
 \lnorm\big( \Phi_n^{-1} A - A_n \Phi_n^{-1}\big)f\lnorm_{\infty}
 \leq \kappa h^\alpha\, \Vert f\Vert^{(2+\alpha)}.
\end{equation}

Now we have a straightforward application of this result on approximations:

\begin{lemma}\label{cor5.1} Let $t\mapsto U(t)$ be a strongly continuous
semigroup in $\dot{C}^{(2+\alpha)}({\bbR}^{d})$, $\Vert U(t)\Vert^{2+\alpha}
\leq \exp(\sigma t)$ with some $\sigma \geq 0$, such that $u(t)=U(t)u_0$ solves
the first IVP in (\ref{exp5.2}). Let $t\mapsto U_n(t)$,
$\lnorm U_n(t)\lnorm_{\infty}\leq 1$, be semigroups generated by $-A_n$ in
$\dot{l}_\infty(G_n)$, such that ${\bf u}_n(t)=U_n(t){\bf u}_{0n}$ solve
the IVPs for ODE in (\ref{exp5.2}). If the differential operator
$A(\msx)$ on ${\bbR}^d$ fulfills the condition (\ref{exp5.5}), then the following
assertion is valid: For each $T > 0$ there exists a positive
number $\rho(T)$ such that
\begin{eqnarray*}
 \sup_{t\in [0,T]}\,\lnorm \Phi_n^{-1}\, u(t) \, -\, {\bf u}_n(t) \lnorm_\infty
 \ \leq \ \rho(T) \| u_0 \|^{2+\alpha} \, h^\alpha,
\end{eqnarray*}
for all $u_0 \in \dot{C}^{(2+\alpha)}({\bbR}^{d})$.
\end{lemma}

{\Proof} The function $s\mapsto {\bf f}(s)=U_n(t-s)\Phi_n^{-1}U(s)u_0\in
\dot{l}_\infty(G_n)$, for $0\leq s\leq t, u_0\in \dot{C}^{(2+\alpha)}({\bbR}^{d})$,
has a continuous derivative of the form ${\bf f}'(s)=U_n(t-s)
\big(A_n \Phi_n^{-1}U(s)-\Phi_n^{-1}AU(s)\big)u_0$. Therefore the following
identity must be valid for each $u_0\in \dot{C}^{(2+\alpha)}({\bbR}^{d})$:
\[ \Big(\Phi_n^{-1} U(t) - U_n(t)\,\Phi_n^{-1}\Big)u_0 \ = \  \int_0^t \ U_n(t-s)
 \, \big (A_n \, \Phi_n^{-1}-\Phi_n^{-1}\, A \big) \, U(s)\,u_0 \, ds.\]
The $\dot{l}_\infty(G_n)$-norm of the integrand is first estimated from above
by
\[\begin{array}{c}
 \lnorm U_n(t-s)\lnorm_\infty\, \lnorm \big(A_n \, \Phi_n^{-1} \,-\, \Phi_n^{-1}
 \, A \big)\, U(s) u_0 \lnorm_\infty \ \leq \ \kappa h^\alpha\,
 \Vert U(s)u_0\Vert^{(2+\alpha)} \\
 \leq \ \kappa h^\alpha\, \exp(\sigma T) \,\Vert u_0\Vert^{(2+\alpha)},\end{array}\]
where (\ref{exp5.5}) is used. Hence,
\[ \blnorm \, \Big( \Phi_n^{-1} U(t) \, - \, U_n(t)\Phi_n^{-1}
 \Big ) u_0\, \blnorm_{\raisebox{-1.4pt}{\mbox{$\infty$}}} \ \leq \
 T\,\exp(\sigma T)\,\kappa\,h^\alpha\,\Vert u_0 \Vert^{2+\alpha},\]
implying the assertion. {\QED}

Due to the density of $\dot{C}^{(2+\alpha)}({\bbR}^{d})$ in $\dot{C}({\bbR}^{d})$
we have the following auxiliary result.

\begin{corollary}\label{cor5.2} Let $A(\msx)$ be defined by (\ref{exp3.1}) and
Assumption \ref{Ass3.1}, and let it fulfil the conditions of Theorem \ref{Th4.1}.
If the coefficients $b_i$ belong to the class $C^{(\alpha)}$ on ${\bbR}^d$
and $a_{ij}$ belong to the class $C^{(1+\alpha)}$ on ${\bbR}^d$ then:
\begin{description}\itemsep -0.1cm
 \item{(i)} The operators $U(n,t)P(n)$ converge strongly in $\dot{C}({\bbR}^{d})$
to $U(t)$, uniformly on segments of $[0,\infty)$.
 \item{(ii)} The limit (\ref{exp1.5}) is valid.
\end{description}
\end{corollary}

The so called classical diffusion, i.e. the process in ${\bbR}^d$ determined by
the differential operator $A(\msx)=-\sum_{ij=1}^da_{ij}(\msx)\partial_i\partial_j$,
is usually simulated by using its representation in trems of stochastic
differential equations. Corollary \ref{cor5.2} makes it possible to simulate
sample paths of classical diffusion in terms of MJPs. This alternative approach
to simulation gives better results in an estimation of the statistical moments of
the first exit time from open sets at subsets of the boundary with a rapidly
changing normal,

\subsection*{Convergence in the general case}
Let us define a mollifier $\msx\mapsto \vartheta(n,\msx) = h^{-d} \vartheta(h^{-1}
\msx)$ in terms of a non-negative function $\vartheta(\cdot)$ of class $C^{(2)}$
on ${\bbR}^d$ with the support equal to the unit ball $B_1(0)\subset {\bbR}^d$
and $\int \vartheta(\msx)
d\msx=1$. A smoothing procedure of coefficients $a_{ij}, b_i$ is determined
by replacing these coefficients with the sequence of coefficients
$a_{ij}^{(m)}=\vartheta(m) \ast a_{ij}, b_i^{(m)}=\vartheta(m) \ast b_i$, $m\in
\bbN$, where $\ast$ denotes the convolution operator.
For each $m\in \bbN$ the resulting differential operator $A^{(m)}(\msx)$
has the closure in $\dot{C}({\bbR}^{d})$ and generates the Feller semigroup
$U^{(m)}(\cdot)$. Now we consider the mappings:
\begin{equation}\label{exp5.7}\begin{array}{ccc}
 A(\msx) &\longmapsto & \{A^{(m)}(\msx)\,:\,m\in \bbN\},\\
 \downarrow && \downarrow\\
 A_n     &\longmapsto & \{A_n^{(m)}\,:\,m\in \bbN\},\end{array}
\end{equation}
and the operators $U^{(m)}(t)$ and $U^{(m)}(n,t)P(n)$ in
the Banach space $\dot{C}({\bbR}^{d})$, where the semigroup $U^{(m)}(\cdot)$ is
generated by the closure of $A^{(m)}(\msx)$, and where $U^{(m)}(n,t)=\exp(-
A_n^{(m)}t)$. For the double sequence of operators
$U^{(m)}(n,t)P(n), m,n\in \bbN$, we will prove the following limits:
\begin{equation}\label{exp5.8}
 U^{(m)}(n,t)P(n) \
 \begin{array}{c} \raisebox{-1.5mm}{$n\to \infty$}\\ \longrightarrow\\  \\
 \end{array} \ U^{(m)}(t) \
 \begin{array}{c} \raisebox{-1.5mm}{$m\to \infty$}\\ \longrightarrow\\  \\
 \end{array} \ U(t),
\end{equation}
uniformly on segments of $[0,\infty)$. By using (\ref{exp5.8})  and applying the
diagonalization argument to the sequence $\{U^{(m)}(n,\cdot)P(m):(m,n)\in
{\bbN}^2\}$, we get the main result of this article:

\begin{theorem}\label{Th5.2} Let the differential operator $A(\msx)$ be
defined by (\ref{exp3.1}) and Assumption \ref{Ass3.1}.
There exists a sequence of pairs $(m,n(m))\in {\bbN}^2$ such that
\[ \lim_m \,\Vert U_{n(m)}^{(m)}(t)P(n(m))-U(t)\Vert_\infty \ = \ 0,\]
uniformly on segments of $[0,\infty)$. Therefore, the asymptotic
(\ref{exp1.5}) is valid for the sequence $\{U_{n(m)}^{(m)}(\cdot):m\in \bbN\}$.
\end{theorem}

{\Proof} Since the coefficients $b_i^{(m)}$ belong to the class $C^{(\alpha)}$ on
${\bbR}^d$ and the coefficients $a_{ij}^{(m)}$ belong to the class $C^{(1+\alpha)}$
on ${\bbR}^d$, the first limit in (\ref{exp5.8}) follows from Corollary \ref{cor5.2}.

It remains to prove the second limit in (\ref{exp5.8}), that is, the
convergence of the semigroups $U^{(m)}(\cdot)$ to $U(\cdot)$ in
$\dot{C}({\bbR}^{d})$, uniformly on segments
of $[0,\infty)$. The convergence of $U^{(m)}(\cdot)$ to $U(\cdot)$
in $L_2({\bbR}^{d})$ can easily be obtained (for instance, by using Theorem 6.1,
Chapter 1 in \cite{EK}). The convergence of $U^{(m)}(t)$ to $U(t)$ in
$\dot{C}({\bbR}^{d})$, uniformly on segments of $[0,\infty)$,
would follow from such convergence on a dense subspace in $\dot{C}({\bbR}^{d})$.
We choose the subspace $C_0^{(\alpha)}({\bbR}^d)=C_0({\bbR}^d)\cap C^{(\alpha)}
({\bbR}^d)$ which is dense in both, the space $L_2({\bbR}^{d})$ and the
space $\dot{C}({\bbR}^{d})$. Due to (i) of Corollary \ref{Cor3.1}, the operators
$U(t)$ and $U^{(m)}(t)$ are bounded in the space $\dot{C}^{(\alpha)}({\bbR}^{d})$,
uniformly on segments $K\subset [0,\infty)$,
i.e. $\Vert U^{(m)}(t)\Vert^{(\alpha)}\leq \beta(K)$, where $\beta(K)$ does not
depend on $m$. Thus we come to the following conclusion. The operators
$U(t)-U^{(m)}(t)$ are continuous mappings from $C_0^{(\alpha)}({\bbR}^d)$ into
$\dot{C}^{(\alpha)}({\bbR}^{d})$, $\Vert U(t)-U^{(m)}(t)\Vert^{(\alpha)}\leq
\beta(K)$, and they converge to zero
in $L_2({\bbR}^{d})$, uniformly on segments $K\subset [0,\infty)$.

Now we apply the following auxiliary result to $u_m(t)=(U(t)-U^{(m)}(t))v,
v\in C_0^{(\alpha)}({\bbR}^d)$. Let $\mathfrak{U} = \{u_n: n \in {\bbN}\}$ be a
sequence of continuous functions on $[0,1]\times {\bbR}^d$ such that:
\begin{description}\itemsep -0.1cm
 \item{a)} $u_n(t)\in L_2({\bbR}^d)$ for each $t\in [0,1]$, and
$\sup\{\lim_n\Vert u_n(t)\Vert_2: t\in [0,1]\}=0$.
 \item{b)} The functions $u_n$ are uniformly H\"{o}lder continuous in the
following sense. There exist $\alpha\in (0,1)$ and $c_\alpha>0$, which do not
depend on $t$ or $\msx$, such that the restrictions $u_n(t)|B_1(\msx)$ fulfil
the following two conditions:
\[ u_n(t)|B_1(\msx) \, \in \, C^{(\alpha)}(B_1(\msx)), \qquad
 \Vert u_n(t) \Vert_\infty^{(\alpha)} \,\leq \,c_\alpha,\]
uniformly with respect to $\msx \in {\bbR}^d$ and $t\in [0,1]$.
\end{description}
Then $\lim_n u_n(t)=0$ in $\dot{C}({\bbR}^d)$, uniformly with respect to $t\in
[0,1]$.

A proof of this auxiliary result is simple.
If the assertion were not valid, there would exist a
positive number $\delta$ and a sequence of pairs $(t_k,\msx_k)\in [0,1]\times
{\bbR}^d, \lim_k|\msx_k|=\infty$, such that $u_k(t_k,\msx_k)\geq \delta$. Due to b)
the following must also be valid: $|u_k(t_k)|\geq\delta/2$ on the ball $B_r(\msx_k)$
where $r =(\delta/2c_\alpha)^{1/\alpha}$. A consequence of these inequalities
would be $\Vert u_k(t_k)\Vert_2 \geq 2^{-1}\delta\sqrt{|B_r({\bf 0})|}$,
contradicting a). {\QED}

We remind the reader that the approximations $A(\msx)\mapsto A^{(m)}(\msx)$,
constructed in the proof of Theorem \ref{Th5.2}, are not the only ones that are
needed for the proof of (\ref{exp1.5}). The approximations defined by
(\ref{exp4.10}) are also needed in order to get matrices $A_n$ of positive type.

\section{SIMULATION OF SAMPLE PATHS}\label{sec6}
Here we demonstrate the efficiency of simulation of sample paths of a
generalized diffusion by using MJPs. We intend to estimate the expectation
and the variance of the first exit time from a Lipshitz domain.
The differential operator $A(\msx)=-\sum_{ij=1}^2 \partial_i a_{ij}(\msx)
\partial_j$ on ${\bbR}^2$ is defined by its diffusion tensor, being a
piecewise constant tensor-valued function of the form,
\[ a(\msx) \ = \ \left [ \begin{array}{cc} \sigma^2 & \alpha(\msx)
 \\ \alpha(\msx) & 1 \end{array} \right ],  \quad \alpha(\msx) \:=\:\rho
 {\bbJ}_{D_0}(\msx), \quad \rho^2 < \sigma^2,\]
where $\sigma^2$ is a positive number, $\rho$ is a real number and $D_0 =
(1/4,3/4)^2$. Let $X(\cdot)$ be diffusion determined by $A(\msx)$ starting
from $\msx_0=(1/2,1/2)$. For a bounded Lipshitz domain $D$ the expectation and
the variance of first exit time from $D$ are given by expressions \cite{Li}:
\begin{equation}\label{exp6.1}
 {\bf E}[\theta] \ = \ \Vert u\Vert_1, \quad
 {\bf Var}[\theta] \ = \ 2\Vert v\Vert_1 \,-\,  {\bf E}[\theta]^2,
\end{equation}
where $u, v$ are the unique solutions of the following boundary value problems,
\begin{equation}\label{exp6.2}
 \begin{array}{cc} \displaystyle
 A(\msx) u(\msx)\ = \ \delta(\msx-\msx_0), & {\bf x} \in D,\\
 u \,|\, \partial D \ = \ 0,& \end{array} \quad
 \begin{array}{cc} \displaystyle
 A(\msx) v(\msx)\ = \ u(\msx), & {\bf x} \in D,\\
 v\,|\, \partial D \ = \ 0,& \end{array}
\end{equation}
and $\delta(\msx)$ is Dirac $\delta$-function at ${\bf 0}$.
We intend to compute ${\bf E}[\theta], {\bf Var}[\theta]$ by simulations
and by using deterministic methods formulated in terms of Expressions
(\ref{exp6.1}), (\ref{exp6.2}). In order to simulate sample paths of $X(\cdot)$
we shall approximate the diffusion by a MJP $X_n(\cdot)$ and simulate
sample paths of $X_n(\cdot)$. In this example we choose $\sigma^2 = 0.1,
\rho =0.02$, $D=(0,1)^2$ and two cases of discretizations, $h=1/200$ and
$h=1/400$.

The generator of the process $X_n(\cdot)$ in $G_n$ is denoted by $A_n({\rm gen})$.
Since the entries $a_{12}$ are non-trivial on $D_0$, we have to use the construction
of Section \ref{sec4} in order to get $A_n({\rm gen})=-A_n$. The parameters
$r_1 = 3, r_2 =1$ of construction are illustrated in Figure~\ref{fig4.2}.
These values of parameters ensure $A_n$ to have the structure of a matrix of
positive type.

Two boundary value problems of (\ref{exp6.2}) have unique solutions
$u, v \in L_1(D)$ as proved in~\cite{BO,LR1}. An efficient numerical method
is constructed and the convergence in $L_1(D)$ is proved in~\cite{LR2}.
This numerical method is based on the construction of $A_n$ which is
described in Section \ref{sec4}.
Efficiency of constructed methods is demonstrated by examples in which
solutions in closed forms are compared with numerical solutions.

Results of computation are expressed in terms of relative errors:
\[ \varepsilon_{exp} \:=\: \frac{<\theta>_{det}-<\theta>_{sim}}{<\theta>_{det}},
 \quad  \varepsilon_{var} \:=\: \frac{\<<\theta\>>_{det}-\<<\theta\>>_{sim}}
 {\<<\theta\>>_{det}},\]
where $<\theta>_{det},<\theta>_{sim}$ are the estimates of ${\bf E}[\theta]$
obtained by using deterministic methods (\ref{exp6.1}) and Monte Carlo simulations,
respectively. Analogously, $\<<\theta\>>_{det},\<<\theta\>>_{sim}$ are the
corresponding quantities for estimates of ${\bf Var}[\theta]$. Some results of
computations are given in the table bellow.
The last column contains the
ratios, $r=t_{det}/t_{sim}$, of computational times $t_{det}$ and $t_{sim}$
of deterministic and Monte Carlo method, respectively. Sample paths are simulated
20000 times.
\vskip0.2cm

\begin{center}
\begin{tabular}{|c|r@{.}l|r@{.}l|r@{.}l|}
\hline
 \rule[-2mm]{0.mm}{6mm}
$h=1/n$& \multicolumn{2}{|c|}{$\varepsilon_{exp}$} &\multicolumn{2}{|c|}
 {$\varepsilon_{var}$} & \multicolumn{2}{|c|}{$r$}\\
\hline
\hline
 n = 200&\rule{0.mm}{4mm} -0&039 & 0&018 & 15&7 \\
 n = 400&                 -0&044 & 0&007 & 20&3 \\
\hline
\end{tabular}
\vskip0.5cm

\centerline{Comparison of results obtained by}
\centerline{deterministic and Monte Carlo methods}
\end{center}
As expected, the first two statistical moments of the first exit time can
be estimated by Monte Carlo simulations dozen times faster than by using
the deterministic method formulated by (\ref{exp6.1}) and (\ref{exp6.2}).

{\bf Acknowledgement}. I would like to thank Vlada Limic for numerous
discussions and helpful comments on previous versions of the manuscript.

\end{document}